\documentclass[reqno]{amsart}
\usepackage{amsfonts,amsmath,amssymb,color}

\numberwithin{equation}{section}

\newtheorem{thm}{Theorem}
\newtheorem{lemma}{Lemma}[section]

\newtheorem{prop}{Proposition}[section]
\newtheorem{defi}{Definition}[section]

\newtheorem{claim}{Claim}[section]

\newcommand{\be}{\begin{equation*}}
\newcommand{\ee}{\end{equation*}}
\newcommand{\beq}{\begin{equation}}
\newcommand{\eeq}{\end{equation}}
\newcommand{\begincal}{\begin{eqnarray*}}
\newcommand{\fincal}{\end{eqnarray*}}
\newcommand{\ds}{\displaystyle}

\newcommand{\ue}{u_\eps}
\newcommand{\R}{{\mathbb R}}
\newcommand{\N}{{\mathbb N}}
\newcommand{\eps}{\varepsilon}
\newcommand{\tue}{\tilde{u}_\varepsilon}
\newcommand{\hue}{\hat{u}_\varepsilon}

\newcommand{\de}{d_\varepsilon}

\newcommand{\rn}{{\mathbb R}^n}
\newcommand{\Om}{\Omega}
\newcommand{\io}{\int_{\Omega}}
\newcommand{\ibo}{\int_{\partial\Omega}}

\begin{document}
\title[Stability of the Poho{\v{z}}aev obstruction]
{Stability of the Poho{\v{z}}aev obstruction\\ in dimension $3$} 

\author{Olivier Druet } 
\address{Olivier Druet, Ecole normale sup\'erieure de Lyon, D\'epartement de 
Math\'ematiques - UMPA, 46 all\'ee d'Italie, 69364 Lyon cedex 07,
France}
\email{Olivier.Druet@umpa.ens-lyon.fr}
%\thanks[Druet]{CNRS - ENS Lyon and PIMS}

\author{Paul Laurain }
\address{Paul Laurain, Ecole normale sup\'erieure de Lyon, D\'epartement de 
Math\'ematiques - UMPA, 46 all\'ee d'Italie, 69364 Lyon cedex 07,
France}
\email{Paul.Laurain@umpa.ens-lyon.fr}
%\thanks[Laurain]{ENS Lyon}

\date{June 2008}

\begin{abstract} We investigate problems connected to the stability of the well-known Poho{\v{z}}aev obstruction. We generalize results which were obtained in the minimizing setting by Brezis and Nirenberg \cite{BrezisNirenberg} and more recently in the radial situation by Brezis and Willem \cite{BrezisWillem}.
\end{abstract}

\maketitle

Let $\Omega$ be a smooth bounded domain in $\rn$, $n\ge 3$. Let $h\in C^1\left(\rn\right)$ and let us consider the equation
\begin{equation}\label{mainequation}
\left\{\begin{array}{l}
{\displaystyle \Delta u + hu = \vert u\vert^{\frac{4}{n-2}} u\hbox{ in }\Omega}\\
\,\\
{\displaystyle u=0\hbox{ on }\partial\Omega}
\end{array}\right.
\end{equation}
where ${\displaystyle \Delta u = -\sum_{i=1}^n \frac{\partial^2 u}{\partial x_i^2}}$. 
It is well-known that, if $\Omega$ is star-shaped with respect to the origin and if $h$ satisfies 
\begin{equation}\label{mainhyp}
h(x)+\frac{1}{2}\left<x, \nabla h(x)\right> \ge 0\hskip.1cm,
\end{equation}
then there are no non-trivial solutions of (\ref{mainequation}). This is a consequence of Poho{\v{z}}aev's identity (see \cite{Pohozaev} and equation (\ref{pohozaev3}) of appendix \ref{Pohozaevidentities}) and is referred to as the Poho{\v{z}}aev obstruction. 

The above equation has been quite intensively studied in the past thirty years. Many existence results have been obtained if $\Omega$ is not assumed to be star-shaped or if $h$ does not verify (\ref{mainhyp}). It is almost impossible to give an exhaustive list of references on this equation.  

\medskip In this paper, we investigate the question of non-existence of positive solutions of equation (\ref{mainequation}) and more precisely the stability properties of the Poho\v{z}aev obstruction. 

\begin{defi}\label{definitionstability}
Let $\Om$ be a star-shaped domain of $\rn$ and let $\left(X, \left\Vert\, .\,\right\Vert_X\right)$ be some Banach space of functions on $\Om$ (typically $X={\mathcal C}^{k,\eta}\left(\Omega\right)$, $X=L^\infty\left(\Om\right)$ or $X=L^p\left(\Omega\right)$). Let $h_0\in X\cap C^1\left(\Om\right)$ be a function which satisfies (\ref{mainhyp}). We say that the Poho{\v{z}}aev obstruction is $X$-stable at $\left(h_0,\Om\right)$ if the following property holds~: there exists $\delta\left(h_0,\Om,X\right)>0$ such that for any function $h\in X$ with 
$$\left\Vert h-h_0\right\Vert_X \le \delta\left(h_0,\Om,X\right)\hskip.1cm,$$
the only non-negative $C^2$-solution of (\ref{mainequation}) is $u\equiv 0$. 

We say that the Poho{\v{z}}aev obstruction is $X$-stable if it is $X$-stable at $\left(h_0,\Om\right)$ for all $\Om$ star-shaped with respect to the origin and all $h_0\in X\cap C^1\left(\Om\right)$ satisfying (\ref{mainhyp}).
\end{defi}

Note that the property (\ref{mainhyp}) is not stable under perturbations of the function $h$ in any $C^k$-space. Since the work of Brezis and Nirenberg \cite{BrezisNirenberg}, we know that equation (\ref{mainequation}) behaves differently in dimension $3$ and in dimensions $n\ge 4$.  It is clear that, in dimensions $n\ge 4$, the Poho{\v{z}}aev obstruction is not $X$-stable for any reasonable $X$. Indeed, any perturbation of $h\equiv 0$ which is negative somewhere leads to a minimizing solution in dimensions $n\ge 4$ (see \cite{BrezisNirenberg})\footnote{Note that this remark concerns only $X$-stability in general. The question of $X$-stability at some given positive function $h$ in dimensions $n\ge 4$ is not investigated in this paper.}. Hence we investigate in this paper the question of the stability of the Poho{\v{z}}aev obstruction for various spaces $X$ in dimension $3$. We give a complete answer to this question in the following theorems. 

\begin{thm}\label{thmstability}
The Poho{\v{z}}aev obstruction is $C^{0,\eta}$-stable for any $\eta>0$ in dimension $3$. In other words, given any $\eta>0$, any domain $\Omega$ in ${\mathbb R}^3$, star-shaped with respect to the origin, and any function $h_0\in C^1\left(\Om\right)$ satisfying (\ref{mainhyp}), there exists $\delta\left(\eta,\Om,h_0\right)>0$ such that, if $h\in C^{0,\eta}\left(\Om\right)$ satisfies 
$$\left\Vert h-h_0\right\Vert_{C^{0,\eta}\left(\Om\right)} \le \delta\left(\eta,\Om,h\right)\hskip.1cm,$$
the only non-negative solution of (\ref{mainequation}) is $u\equiv 0$. 
\end{thm}

Note that a consequence of our theorem is the following~: if $\Om$ is a star-shaped domain in ${\mathbb R}^3$, there exists a constant $\hat{\lambda}\left(\Omega\right)>0$ such that equation (\ref{mainequation}) does not possess any nontrivial positive solutions with $h\equiv \lambda$ for $\lambda>-\hat{\lambda}\left(\Omega\right)$. This is in sharp contrast with the situation for non star-shaped domains (see \cite{BahriCoron} for instance). 

In the seminal paper \cite{BrezisNirenberg}, it was proved that there are no minimizing solutions of equation (\ref{mainequation}) in dimension $3$ if the function $h\ge -\lambda^\star\left(\Omega\right)$ for some $\lambda^\star\left(\Omega\right)>0$. Since $h\ge 0$ if $h$ satisfies (\ref{mainhyp}), a  consequence of this result is a version of the above stability in $C^0$ when one considers only minimizing solutions. A necessary and sufficient condition on the function $h$ and the domain $\Omega$ for the existence of a minimizing solution of (\ref{mainequation}) in dimension $3$ was found in \cite{Druet}. 

In \cite{BrezisWillem}, the authors studied this question in the case of the unit ball with radial functions. If we let 
$$L^p_r \left(B\right)= \left\{u\in L^p\left(B\right)\,,\,\, u\hbox{ radial }\right\}\hskip.1cm,$$
then it was proved in \cite{BrezisWillem} that the Poho{\v{z}}aev obstruction is $L^\infty_r$-stable\footnote{One should restrict oneself to radial solutions of the equation in the definition of stability.} on the unit ball of ${\mathbb R}^3$ for all functions $h\in L^\infty_r\left(B\right)\cap C^1\left(B\right)$. In \cite{BrezisWillem}, the question of the extension of the result to the non-radial case was explicitly asked. Our result provides an answer to this question. However, the situation is more delicate than expected in the non-radial case since, while the Poho{\v{z}}aev obstruction is $C^{0,\eta}$-stable for all $\eta>0$, it is never $L^\infty$-stable.

\begin{thm}\label{thminstability}
The Poho{\v{z}}aev obstruction is never $L^\infty$-stable. In other words, given any $\eps>0$, any domain $\Omega$ in ${\mathbb R}^3$, star-shaped with respect to the origin and any function $h_0\in C^1\left(\Om\right)$ satisfying (\ref{mainhyp}), we can find some function $h_\eps\in L^\infty\left(\Omega\right)$ such that 
$$\left\Vert h_\eps-h\right\Vert_\infty\le \eps$$
and some positive functions $u_\eps\in C^2\left(\Omega\right)$ satisfying the equation 
$$\left\{\begin{array}{l}
{\displaystyle \Delta u_\eps + h_\eps u_\eps = u_\eps^5\hbox{ in }\Omega}\\
\,\\
{\displaystyle u_\eps=0\hbox{ on }\partial\Omega\, ,\,\, u_\eps> 0\hbox{ in }\Omega}
\end{array}\right.$$
\end{thm}

\medskip Thus the $L^\infty_r$-stability result obtained by Brezis-Willem is really specific to the radial case. In fact, it is not really due to the symmetry of the solutions but to one of its by-product in dimension $3$, precisely that sequences of solutions of equation (\ref{mainequation}) which are radial are either compact or develop only one concentration point. In fact, with the PDE techniques (to be compared to the ODE techniques used in \cite{BrezisWillem}) we use below, we can revisit the question of the stability of the Poho{\v{z}}aev obstruction in dimension $3$ in the radial case. We improve the result of \cite{BrezisWillem} by proving that the Poho\v{z}aev obstruction is $L^p_r$-stable on the unit ball for all $p>3$ but is never $L^3_r$-stable. For precise statements, we refer the reader to the end of section \ref{sectionstability} and the beginning of section \ref{sectionexamples}.

\medskip All these results give a complete picture of the stability of the Poho{\v{z}}aev obstruction in dimension $3$ when the attention is restricted to non-negative solutions. The question remains widely open if one allows solutions to change sign, and is certainly more subtle due to the variety of changing-sign solutions of $\Delta u = u^5$ in ${\mathbb R}^3$. 

\medskip The paper is organized as follows. Section \ref{sectionstability} is devoted to the proofs of theorem \ref{thmstability} and of the corresponding result in the radial situation. The proof makes use of standard blow-up analysis in dimension $3$ (see section \ref{sectionblowupanalysis}) and of an extension of Poho{\v{z}}aev's identity to Green's functions (see Appendix \ref{pohog}). Section \ref{sectionexamples} is devoted to the proofs of theorem \ref{thminstability}  and of the corresponding result in the radial situation. Here we have to construct examples of functions $h$ arbitrarily close in $X$ to some given function for which there is a positive solution of equation (\ref{mainequation}). It appears to be quite subtle because we need to be sharp. For instance, in order to prove theorem \ref{thminstability}, our functions $h$ must be close to the given function in $L^\infty\left(\Om\right)$ but not in $C^{0,\eta}\left(\Om\right)$ for any $\eta>0$. 

\medskip {\bf Acknowledgements} : The first author wishes to thank H. Brezis for putting to his attention the question raised in \cite{BrezisWillem}. This work was done while the two authors were staying at PIMS, Vancouver, for one year. They wish to thank all the members of PIMS for their great hospitality during this year.

\section{Pointwise analysis around a concentration point}\label{sectionblowupanalysis}

We consider in this section a sequence $\left(h_\eps\right)$ in $C^{0,\eta}\left(\rn\right)$ for some $\eta>0$ and a sequence $\left(u_\eps\right)$ of $C^2$-solutions of 
\begin{equation}\label{eq1.1}
\left\{\begin{array}{ll}
{\ds \Delta u_\eps + h_\eps u_\eps =u_\eps ^5}&{\ds \hbox{ in } \Omega}\\
\, &\,\\
{\ds u_\eps =0 }&{\ds \hbox{ on } \partial\Omega}\\
\, &\, \\
{\ds \ue >0}&{\ds \hbox{ in } \Omega}
\end{array}\right.
\end{equation}
where $\Omega$ is some smooth domain of $\R^3$ and
\begin{equation}\label{eq1.2}
h_\eps \rightarrow h \hbox{ in } L^{p}(\Omega) \hbox{ as } \eps \rightarrow 0
\end{equation}
for some $p>3$ where $h\in C^1(\R^3)$ satisfies $h\geq 0$ in $\Omega$. Note that, as soon as $h$ satisfies (\ref{mainhyp}), it is non-negative. 

We also assume that we have a sequence $\left(x_\eps\right)$ of points in $\Omega$ and a sequence $\left(\rho_\eps\right)$ of positive real numbers with $0< 3 \rho_\eps \leq d(x_\eps, \partial\Omega)$  such that 
\beq\label{eq1.3}
\nabla \ue( x_\eps) =0
\eeq
and
\beq\label{eq1.4}
 \rho_\eps \left[ \sup_{B(x_\eps, \rho_\eps)}  u_\eps(x) \right]^2  \rightarrow + \infty \hbox{ as }\eps \rightarrow 0 \hskip.1cm.
\eeq
We prove in this section that the following holds~:

\begin{prop}\label{estim}
If there exists $C_0 >0$ such that
\beq\label{eq1.5}
\vert x_\eps - x \vert^\frac{1}{2} u_\eps\leq C_0 \hbox{ in } B(x_\eps, 3\rho_\eps)\hskip.1cm,
\eeq
then there exists $C_1 >0$ such that
\be
\begin{split}
&u_\eps(x_\eps) u_\eps(x) \leq C_1 \vert x_\eps - x \vert^{-1}\hbox{ in } B(x_\eps, 2\rho_\eps)\setminus \{x_
\eps \}\hbox{ and }\\
&u_\eps(x_\eps) \vert\nabla u_\eps(x)\vert \leq C_1 \vert x_\eps - x \vert^{-2}\hbox{ in } B(x_\eps, 2\rho_\eps)\setminus \{x_
\eps \}.
\end{split}
\ee
Moreover, if $\rho_\eps \rightarrow 0$, then 
\be 
\rho_\eps u_\eps(x_\eps) u_\eps(x_\eps+\rho_\eps x) \rightarrow \frac{1}{\vert x\vert}+ b  \hbox{ in } C^1_{loc}(B(0, 2)\setminus \{0\})\hbox{ as }\eps\to 0
\ee
where $b$ is some harmonic function in $B(0,2)$ with $b(0)=0$. At last, if the convergence in (\ref{eq1.2}) holds in $C^{0,\eta}$, then we also get that $\nabla b(0)=0$.
\end{prop}

\medskip The rest of this section is dedicated to the proof of this proposition. We follow the lines of \cite{DruetIMRN}, section 2 (see also \cite{DruetHebey08}). However, one must note that, compared to \cite{DruetHebey08} and other works on this kind of blow-up analysis, some new difficulties arise since  the linear term $\left(h_\eps\right)$ is only uniformly bounded in some $L^p\left(\Omega\right)$.

We divide the proof of the proposition into several claims. The first one gives the asymptotic behaviour of $u_\eps$ around $x_\eps$ at an appropriate small scale. 

\begin{claim}
\label{estim1}
After passing to a subsequence, we have that
\beq\label{eq1.7}
\mu_\eps^\frac{1}{2}  u_\eps(x_\eps + \mu_\eps x) \rightarrow \frac{1} { \left(1+ \frac{\vert x \vert^{2} }{3} \right)^{\frac{1}{2}} }
\hbox{ in } C^1_{loc}(\R^3), \hbox{ as } \eps \rightarrow 0
\eeq
where $\mu_\eps = u_\eps\left(x_\eps\right)^{-2}$.
\end{claim}

\medskip {\bf Proof of Claim \ref{estim1}.} Let $\tilde{x}_\eps \in \overline{B(x_\eps, \rho_\eps)}$ and $\tilde{\mu}_\eps>0$ be such that
\beq\label{eq1.8}
u_\eps(\tilde{x}_\eps) = \sup_{B(x_\eps, \rho_\eps)} u_\eps=\tilde{\mu}_\eps^{-\frac{1}{2}} \hskip.1cm.
\eeq
Thanks to (\ref{eq1.4}), we have that 
\beq\label{eq1.9}
\tilde{\mu}_\eps\to 0 \hbox{ and }\frac{\rho_\eps}{\tilde{\mu}_\eps}\to +\infty\hbox{ as }\eps\to 0\hskip.1cm.
\eeq
Thanks to  (\ref{eq1.5}), we also have that
\beq
\label{eq1.10}
\vert x_\eps-\tilde{x}_\eps\vert= O(\tilde{\mu}_\eps).
\eeq
We set for $\displaystyle x \in \Omega_\eps =\left\{ x\in\R^3 \hbox{ s.t. } \tilde{x}_\eps + \tilde{\mu}_\eps x \in \Omega \right\}$, 
\be
\tue (x)= \tilde{\mu}_\eps^{\frac{1}{2}}  u_\eps(\tilde{x}_\eps + \tilde{\mu}_\eps x) 
\ee
which verifies
\beq
\label{eq1.11}
\begin{split}
&\Delta \tue +\tilde{\mu}_\eps^2 \tilde{h}_\eps \tue =\tue ^5 \hbox{ in } \Omega_\eps\hskip.1cm, \\
&\tue(0) = \sup_{ B(\frac{x_\eps -\tilde{x}_\eps}{\tilde{\mu}_\eps}, \frac{\rho_\eps}{\tilde{\mu}_\eps})} \tue =1\hskip.1cm,
\end{split}
\eeq
where $\tilde{h}_\eps= h\left(\tilde{x}_\eps + \tilde{\mu}_\eps x\right)$. Thanks to (\ref{eq1.4}), (\ref{eq1.8}) and (\ref{eq1.10}), we get that
\beq
\label{eq1.12}
B\left(\frac{x_\eps -\tilde{x}_\eps}{\tilde{\mu}_\eps}, \frac{\rho_\eps}{\tilde{\mu}_\eps}\right) \rightarrow \R^3 \hbox{ as } \eps \rightarrow 0 \hskip.1cm.
\eeq
Now, thanks to (\ref{eq1.11}), (\ref{eq1.12}), and by standard elliptic theory, we get that, after passing to a subsequence, $\tue\rightarrow U$ in $C^1_{loc}(\R^3)$ as $\eps \rightarrow 0$ where $U$ satisfies
\be 
\Delta U = U^5 \hbox{ in } \R^3 \hbox{ and } 0\leq U \leq1=U(0) \hskip.1cm.
\ee
Thanks to the work of Caffarelli, Gidas and Spruck \cite{CGS}, we know that
\be
U(x) = \left(1+ \frac{\vert x \vert^{2} }{3} \right)^{-\frac{1}{2}} \hskip.1cm.
\ee
Moreover, thanks to (\ref{eq1.10}), we know that, after passing to a new subsequence, $\frac{x_\eps -\tilde{x}_\eps}{\tilde{\mu}_\eps} \rightarrow x_0 \hbox{ as } \eps \rightarrow 0$ for some $x_0 \in \R^3$. Hence, since $x_\eps$ is a critical point of $u_\eps$, $x_0$ must be a critical point  of $U$, namely  $x_0=0$. We deduce that $\frac{\mu_\eps}{\tilde{\mu}_\eps} \rightarrow 1$ where $\mu_\eps$ is as in the statement of the claim. The claim \ref{estim1} follows. \hfill $\blacksquare$

\medskip For $0\leq r \leq 3 \rho_\eps$, we set
$$\psi_\eps(r)= \frac{r^{\frac{1}{2}} }{\omega_2 r^2} \int_{\partial B({x_\eps},r)} u_\eps d\sigma$$
where $d\sigma$ denotes the Lebesgue measure on the sphere $\partial B({x_\eps},r)$ and $\omega_2=4\pi$ is the volume of the unit $2$-sphere. 
We easily check, thanks to Claim \ref{estim1}, that
\beq\label{eq1.13}
\psi_\eps(\mu_\eps r)= \left( \frac{r}{1+ \frac{r^2}{3}}\right)^\frac{1}{2} + o(1)\hbox{, }\psi_\eps ' (\mu_\eps r)= \frac{1}{2} \left( \frac{r}{1+ \frac{r^2}{3}}\right)^\frac{3}{2} \left(\frac{1}{r^2} - \frac{1}{3} \right) + o(1)\hskip.1cm.
\eeq
We define $r_\eps$ by
\be
r_\eps = \max \left\{ r\in[ 2\sqrt{3}\mu_\eps, \rho_\eps] \hbox{ s.t. }   \psi_\eps ' (s) \leq  0 \hbox{ for } s \in [2\sqrt{3} \mu_\eps, r] \right\} \hskip.1cm.
\ee
Thanks to (\ref{eq1.13}), the set on which the maximum is taken is not empty for $\eps$ small enough, and moreover 
\beq\label{eq1.14}
\frac{r_\eps}{\mu_\eps} \rightarrow + \infty   \hbox{ as } \eps \rightarrow 0\hskip.1cm.
\eeq
We prove now the following~: 

\begin{claim} \label{estim2} 
There exists $C >0$, independent of $\eps$, such that
\be
\begin{split}
&u_\eps(x) \leq C \mu_\eps^\frac{1}{2} \vert x_\eps - x \vert^{-1}\hbox{ in } B(x_\eps, 2r_\eps)\setminus \{x_
\eps \}\hbox{ and } \\
&\vert \nabla u_\eps(x) \vert \leq C  \mu_\eps^\frac{1}{2} \vert x_\eps - x \vert^{-2}\hbox{ in } B(x_\eps, 2r_\eps)\setminus \{x_
\eps \}\hskip.1cm.
\end{split}
\ee
\end{claim}

\medskip {\bf Proof of Claim \ref{estim2}.} We follow the proof of Lemma 1.5 and 1.6 of \cite{DruetHebey08}. However, there is an extra-difficulty due to the fact that we do not assume any pointwise convergence of $h_\eps$ to $h$. We first prove that for any given $0< \nu <\frac{1}{2}$, there exists $C_\nu >0 $ such that 
\beq\label{eq1.15}
\ue(x) \leq C_\nu \left( \mu_{\eps}^{\frac{1}{2}(1-2\nu)} \vert x-x_\eps \vert^{-(1-\nu)} + \alpha_\eps 
\left(\frac{r_{\eps}}{\vert x-x_\eps \vert}\right)^{\nu} \right)
\eeq
for all $x\in B\left(x_{\eps},2r_{\eps}\right)$ and $\eps$ small enough, where
\beq\label{eq1.15bis}
\alpha_\eps= \left( \sup_{ \partial B\left(x_\eps, r_\eps\right)} \ue\right)\hskip.1cm.
\end{equation}
First of all, we can use (\ref{eq1.5}) and apply the Harnack inequality, see for instance theorem 4.17 of \cite{HanLin}, to get the existence of some $C>0$ such that 
\beq\label{eq1.16}
\frac{1}{C} \max_{\partial B(x_\eps, r)} \left(\ue+r\left\vert \nabla u_\eps\right\vert\right) \leq \frac{1}{\omega_2 r^2} \int_{\partial B({x_\eps},r)} u_\eps d\sigma \leq C \min_{\partial B(x_\eps, r)}\ue  
\eeq
for all  $0  <r < \frac{5}{2} \rho_\eps $ and all $ \eps>0$. The details of the proof of such an assertion may be found in \cite{DruetHebey08}, lemma 1.3. Hence, thanks to (\ref{eq1.13}) and (\ref{eq1.14}), we have that 
$$\vert x-x_\eps\vert^\frac{1}{2} \ue (x) \leq C \psi_\eps(r) 
\leq C \psi_\eps(R \mu_{\eps})=C \left( \frac{R}{1+ \frac{R^2}{3}}\right)^\frac{1}{2} + o(1)$$
for all $R\geq  2\sqrt{3}$, all $r\in[ R \mu_{\eps},r_\eps ]$, all $\eps$ small enough and all $x\in \partial B\left(x_{\eps},r\right)$. Thus we get that 
\beq\label{eq1.17}
\sup_{B(x_\eps, r_\eps)\setminus B(x_\eps, R \mu_\eps)}  \vert x-x_\eps\vert^\frac{1}{2} \ue (x) = e(R)+o(1) 
\eeq
where $e(R)\rightarrow 0$ as $R \rightarrow +\infty$. Let $0<\sigma\le 1$ and ${\mathcal G}_{\eps,\sigma}$ be the Green function of the operator $\Delta + \frac{ h_{\eps}}{\sigma}$ in $\Omega$ with Dirichlet boundary condition. Thanks to the fact that $h$ is non-negative (this is an assumption in this section), we can use lemma \ref{estimg} to get the existence of some $C_\sigma >0$ such that 
\beq\label{eq1.18}
\left\vert \vert x-y\vert {\mathcal G}_{\eps, \sigma}(x,y) -\frac{1}{\omega_2} \right\vert \leq C_{\sigma} \vert x-y \vert ,
\eeq 
and that
 \beq
 \label{eq1.19}
 \left\vert \vert x-y\vert^2 \vert\nabla {\mathcal G}_{\eps, \sigma}(x,y)\vert -\frac{1}{\omega_2} \right\vert \leq C_{\sigma} \vert x-y \vert ,
 \eeq 
 for all $x\not=y \in \Omega$. We fix $0<\nu < \frac{1}{2}$ and we set
 \be 
 \Phi_{\eps,\nu} = \mu_\eps^{\frac{1}{2}(1-2\nu)} {\mathcal G}_{\eps, 1-\nu}(x_\eps,x)^{1-\nu}+\alpha_\eps \bigl(r_\eps {\mathcal G}_{\eps, \nu}(x_\eps,x)\bigr)^{\nu} .
 \ee
Thanks to (\ref{eq1.18}), (\ref{eq1.15})  reduces to prove that
\be
\sup_{B(x_\eps, 2 r_\eps)} \frac{\ue}{ \Phi_{\eps,\nu}} = O(1)\hskip.1cm.
\ee
We let $y_\eps \in \overline{B(x_\eps, 2 r_\eps)\setminus\{x_{\eps}\}} $ be such that  
\be
\sup_{B(x_\eps, 2 r_\eps)} \frac{\ue}{ \Phi_{\eps,\nu}} =  \frac{\ue(y_\eps)}{ \Phi_{\eps,\nu}(y_\eps)}.
\ee
We are going to consider the several possible beahviour of the sequence $\left(y_\eps\right)$.

First of all, assume that
\be
\frac{\vert x_{\eps} -y_\eps\vert}{\mu_\eps} \rightarrow R \hbox{ as } \eps \rightarrow 0\hskip.1cm.
\ee
Thanks to Claim \ref{estim1}, we have in this case that  
\be
\mu_\eps^\frac{1}{2}  \ue(y_\eps) \rightarrow (1+R^2)^{-\frac{1}{2}} \hbox{ as } \eps \rightarrow 0.
\ee
On the other hand, thanks to (\ref{eq1.17}), we can write that 
\begincal
\mu_\eps^\frac{1}{2}  \Phi_{\eps,\nu}(y_\eps) &= &\left( \frac{\mu_\eps}{\omega_2 \vert x_{\eps} -y_\eps\vert}\right)^{1-\nu} +O\left(\alpha_\eps \mu_\eps^\frac{1}{2} \left(\frac{r_\eps}{\vert x_{\eps} -y_\eps\vert}\right)^\nu\right) + o(1) \\
 &=&\left( R\omega_2\right)^{\nu-1} + O\left((r_\eps^{\frac{1}{2}} \alpha_\eps) \mu_\eps^{\frac{1}{2}(1-2\nu)} r_\eps^{\frac{1}{2}(2\nu-1)} \right)+ o(1)  \\
 &=&\left( R\omega_2\right)^{\nu-1} + o(1) ,
\fincal
if $R>0$, and $\mu_\eps^\frac{1}{2}  \Phi_{\eps,\nu}(y_\eps) \rightarrow +\infty$ as $\eps \rightarrow 0$ if $R=0$. In any case, $\left(\frac{\ue(y_\eps)}{ \Phi_{\eps,\nu(y_\eps)}}\right)$ is bounded. 

Assume now that there exists $\delta>0$ such that $y_\eps \in B(y_\eps ,r_\eps)\setminus B(y_\eps ,\delta r_\eps)$. Thanks to Harnack's inequality (\ref{eq1.16}), we get that $\ue(y_\eps) =O(\alpha_\eps)$ which, thanks to (\ref{eq1.18}), easily gives that $\frac{\ue(y_\eps)}{ \Phi_{\eps,\nu(y_\eps)}} = O \left(1\right)$.

Hence, we are left with the following situation~:
\beq\label{eq1.20}
 \frac{\vert y_\eps -x_\eps\vert}{r_\eps} \rightarrow 0  \hbox{ and } \frac{\vert x_{\eps} -y_\eps\vert}{\mu_\eps} \rightarrow +\infty \hbox{ as } \eps \rightarrow 0\hskip.1cm.
\eeq
Thanks to the definition of $y_\eps$, we can then write that
\be
\frac{\Delta u_\eps(y_\eps)}{u_\eps(y_\eps)} \geq \frac{\Delta \Phi_{\eps, \nu}(y_\eps)}{\Phi_{\eps, \eta} (y_\eps)}
\ee
which gives, thanks to the definition of $\Phi_{\eps, \nu}$ and multiplying by $\vert x_\eps -y_\eps \vert^2$, that
\begincal
\vert x_\eps -y_\eps \vert^2 u_\eps(y_\eps)^4 &\geq & \nu(1-\nu)\frac{\vert x_\eps -y_\eps \vert^2 }{\Phi_{\eps, \eta} (y_\eps)} \left(  \alpha_\eps r_\eps^\nu \frac{\vert \nabla G_{\eps,\nu}(x_\eps,y_\eps)\vert^2}{G_{\eps, \nu}(x_\eps,y_\eps)^2} G_{\eps,\nu}(x_\eps,y_\eps)^{\nu} \right. \\
&&\quad  +\left. \mu_\eps^{\frac{1}{2}(1-2\eta)} \frac{\vert \nabla G_{\eps, 1-\nu}(x_\eps,y_\eps)\vert^2}{G_{\eps, 1-\nu}(x_\eps,y_\eps)^2} G_{\eps, 1-\nu}(x_\eps,y_\eps)^{1-\nu}\right)\hskip.1cm .
\fincal
Here is the main difference with \cite{DruetHebey08}. Thanks to our choice of  $\Phi_{\eps, \nu}$, the terms involving $h_\eps$ disappear, which is necessary since we did not assume any pointwise convergence of $h_\eps$. Thanks to (\ref{eq1.17}), the left-hand side goes to $0$ as $\eps \rightarrow 0$. Then, thanks to (\ref{eq1.18}), (\ref{eq1.19}) and (\ref{eq1.20}), we get that 
\be
o(1)\geq \nu(1-\nu) + o(1)
\ee
which is a contradiction, and shows that this last case can not occur. This ends the proof of (\ref{eq1.15}).

\medskip We now claim that there exists $C >0$, independent of $\eps$, such that 
\beq
\label{eq1.21}
\ue(x) \leq C \left( \mu_{\eps}^{\frac{1}{2}} \vert x-x_\eps \vert^{-1} +  \alpha_\eps \right) \hbox{ in }  B(x_\eps, r_\eps)\hskip.1cm.
\eeq
Thanks to Claim \ref{estim1} and (\ref{eq1.16}), this holds for all sequences $y_\eps \in B(x_\eps, r_\eps)\setminus\{ x_\eps \} $ such that $\vert y_\eps -x_\eps\vert=O(\mu_\eps)$ or $\frac{\vert y_\eps -x_\eps\vert}{r_\eps} \not\rightarrow 0$. Thus we may assume from now that
\be
\frac{\vert y_\eps -x_\eps\vert}{\mu_\eps} \rightarrow +\infty \hbox{ and }\frac{\vert y_\eps -x_\eps\vert}{r_\eps} \rightarrow 0 \hbox{ as } \eps \rightarrow 0 \hskip.1cm.
\ee
Thanks to the Green representation formula, we write with (\ref{eq1.18}) and (\ref{eq1.19}) that 
\begincal 
\ue(y_\eps)&=& \int_{B(x_\eps, r_\eps)} \mathcal{G}_{\eps,1} (\Delta \ue +h_\eps \ue) \,dx \\
&&+ O\left(r_\eps^{-1}\int_{\partial B(x_\eps, r_\eps)} \left\vert {\partial_\nu} \ue\right\vert \,d\sigma\right) \\
&&+ O\left(r_\eps^{-2}\int_{\partial B(x_\eps, r_\eps)} u_\eps \,d\sigma\right)\hskip.1cm.
\fincal
This gives with (\ref{eq1.15bis}), (\ref{eq1.16}) and (\ref{eq1.18}) that 
\beq
\label{eq1.22}
\ue(y_\eps)= O\left(\int_{B(x_\eps, r_\eps)} \vert x-y_\eps \vert^{-1} \vert \Delta \ue + h_\eps \ue \vert dx\right) + O \left( \alpha_\eps\right) \hskip.1cm.
\eeq  
Using (\ref{eq1.15}) with $\nu=\frac{1}{5}$, we can write that 
\be
\begin{split}
&\int_{B(x_\eps, r_\eps)} \vert x-y_\eps \vert^{-1} \vert \Delta \ue +h_\eps \ue \vert dx \\
&\quad =
\int_{B(x_\eps, \mu_\eps)} \frac{\ue^5}{\vert x-y_\eps \vert}  dx  + \int_{B(x_\eps, r_\eps)\setminus B(x_\eps, \mu_\eps)} \vert x-y_\eps \vert^{-1} \ue^5  dx  \\
& \quad= O\left(\mu_{\eps}^{\frac{1}{2}} \vert y_\eps-x_\eps \vert^{-1} \right)  +  \alpha_\eps^5 r_\eps  \int_{B(x_\eps, r_\eps)\setminus B(x_\eps, \mu_\eps)} \vert x-y_\eps \vert^{-1}   \vert x-x_\eps \vert^{-1}\, dx\\
&\qquad + \mu_{\eps}^{\frac{3}{2}}\int_{B(x_\eps, r_\eps)\setminus B(x_\eps, \mu_\eps)} \vert x-y_\eps \vert^{-1}   \vert x-x_\eps \vert^{-4}dx  \\
&\qquad  \\
& \quad = O\left(\mu_{\eps}^{\frac{1}{2}} \vert y_\eps-x_\eps \vert^{-1} \right) +O( \alpha_\eps^5 r_\eps^{2})\hskip.1cm.
\end{split}
\ee
Thanks to (\ref{eq1.14}) and to (\ref{eq1.17}), this leads to  
\be
\int_{B(x_\eps, r_\eps)} \vert x-y_\eps \vert^{-1} \vert \Delta \ue \vert dx \leq O(\mu_{\eps}^{\frac{1}{2}} \vert y_\eps-x_\eps \vert^{-1} ) +o( \alpha_\eps),
\ee
which, thanks to (\ref{eq1.22}), proves (\ref{eq1.21}). 

In order to end the proof of the first part of the claim, we just have to prove that
\beq
\label{eq1.23}
\alpha_\eps =\sup_{ \partial B(x_\eps, r_\eps)} \ue = O \left( \mu_{\eps}^{\frac{1}{2}} r_\eps^{-1}  \right)\hskip.1cm.
 \eeq
For that purpose, we use the definition of $r_\eps$ to write that 
\be 
(\beta r_\eps)^\frac{1}{2} \psi_\eps(\beta r_\eps) \geq  ( r_\eps)^\frac{1}{2} \psi_\eps( r_\eps) 
\ee
for all $0<\beta <1$. Using (\ref{eq1.16}), this leads to
\be 
r_\eps^\frac{1}{2} \left(  \sup_{ \partial B(x_\eps, r_\eps)} \ue \right) \leq C (\beta r_\eps)^\frac{1}{2} \left(  \sup_{ \partial B(x_\eps, \beta r_\eps)} \ue \right) \hskip.1cm.
\ee 
Thanks to (\ref{eq1.21}), we obtain that 
\be 
\left(  \sup_{ \partial B(x_\eps, r_\eps)} \ue \right) \leq C \beta ^\frac{1}{2} \left(\mu_{\eps}^\frac{1}{2} (\beta r_\eps)^{-1}  +\sup_{ \partial B(x_\eps,  r_\eps)} \ue \right) \hskip.1cm.
\ee 
Choosing $\beta$ small enough clearly gives (\ref{eq1.23}) and thus the pointwise estimate on $\ue$ of the claim. The estimate on $\nabla \ue$ then follows from standard elliptic theory.\hfill$\blacksquare$

\medskip We now prove the following~:

\begin{claim}\label{estim3}
If $r_\eps \rightarrow 0$ as $\eps\rightarrow 0$, up to passing to a subsequence,
\be
r_\eps u_\eps(x_\eps) u_\eps(x_\eps+r_\eps x) \rightarrow \frac{1}{ \vert x\vert}+ b \hbox{ in } C_{loc}^1\left(B\left(0, 2\right)\setminus\{0\}\right)\hbox{ as } \eps \rightarrow 0
\ee
where $b$ is some harmonic function in $B(0,2)$. Moreover, if $r_\eps < \rho_\eps$, then 
$b(0)=1$.
\end{claim}

\medskip {\bf Proof of Claim \ref{estim3}.} We set, for $x \in B(0,2)$,
\be
\tue (x)= \mu_\eps^{-\frac{1}{2}}  r_\eps u_\eps(x_\eps + r_\eps x) 
\ee
which verifies
\beq
\label{eq1.24}
\Delta \tue +r_\eps^2 \tilde{h}_\eps \tue =\left( \frac{\mu_\eps}{r_\eps}\right)^2 \tue ^5 \hbox{ in }  B(0,2) 
\eeq
 where $\tilde{h}_\eps= h(x_\eps + r_\eps x)$. Thanks to Claim \ref{estim2}, there exists $C>0$ such that
 \beq
 \label{eq1.25}
 \tue(x) \leq   \frac{C}{\vert x\vert} \hbox{ in } B(0,2)\setminus \{0\} .
\eeq
Then, thanks to standard elliptic theory, we get that, after passing to a subsequence, $\tue\rightarrow U$ in $C^1_{loc}\left(B(0,2)\setminus \left\{0\right\}\right)$ as $\eps \rightarrow 0$ where $U$ is a non-negative solution of
\be
\Delta U = 0 \hbox{ in } B(0,2) \setminus \{0\} \hskip.1cm.
\ee
Then, thanks to the B\^{o}cher theorem on singularities of harmonic functions, we get that  
\be
U(x)=\frac{\lambda}{\vert x\vert} +b(x)
\ee
where $b$ is some harmonic function in $B(0,2)$ and $\lambda\ge 0$. Now, integrating (\ref{eq1.24}) on $B\left(0,1\right)$, we get that 
\be
\int_{\partial B(0,1)}  \partial_\nu \tue d\sigma =\int_{B(0,1)}  \left(r_\eps^2 \tilde{h}_\eps \tue -\left( \frac{\mu_\eps}{r_\eps}\right)^2 \tue^5 \right)dx
\ee
Thanks to Claim \ref{estim2}, 
\be
\int_{B(0,1)}  r_\eps^2 \tilde{h}_\eps \tue dx \rightarrow 0 \hbox{ as } \eps \rightarrow 0 
\ee
and, thanks to Claim \ref{estim1}, 
\be
\int_{B(0,1)}  \left( \frac{\mu_\eps}{r_\eps}\right)^2 \tue^5 dx \rightarrow \int_{\R^3}  \left(1+\frac{\vert x\vert^2}{3} \right) ^{-\frac{5}{2}} dx = \omega_2 \hbox{ as } \eps \rightarrow 0\hskip.1cm.
\ee
On the other hand, we have that 
\be
\int_{\partial B(0,1)}  \partial_\nu \tue d\sigma \rightarrow - \omega_2 \lambda \hbox{ as } \eps \rightarrow 0\hskip.1cm.
\ee
We deduce that $\lambda=1$, which proves the first part of the Claim.

\medskip Now, if $r_\eps <\rho_\eps$, we have thanks to the definition of $r_\eps$ that
\be
\psi_\eps'(r_\eps) =0\hskip.1cm.
\ee
Setting  $\tilde{\psi}_\eps(r) = \left( \frac{r_\eps}{\mu_\eps}\right)^{\frac{1}{2}} \psi_\eps(r_\eps r)$ for $0<r<2$, we see that
\be
\tilde{\psi}_\eps(r) \rightarrow  \frac{r^\frac{1}{2}}{\omega_2 r^2} \int_{\partial B(0,r)} U d\sigma =r^{-\frac{1}{2}} + r^{\frac{1}{2}} b(0)\hskip.1cm. 
\ee
We deduce that $b(0)=1$, which ends the proof of the Claim.\hfill $\blacksquare$

\medskip We prove at last the following~:

\begin{claim}
\label{estim4}
Using the notations of Claim \ref{estim3}, we have that $b(0)=0$, and if the convergence in (\ref{eq1.2}) holds in $C^{0,\eta}$, then $\nabla b (0) =0$.
\end{claim}

\medskip {\bf Proof of Claim \ref{estim4}.} We use the notation of the proof of  Claim \ref{estim3}.
Let us apply the Poho\v{z}aev identity (\ref{pohozaev1}) of appendix \ref{Pohozaevidentities} to $\tue$ in $B(0,1)$. We obtain that 
\be
 \frac{1}{2} \int_{B(0,1)} r_\eps^2 \left( \tilde{h}_\eps \tue^2 + \tilde{h}_\eps < x ,\nabla \tue^2 >\right) dx = 
 \tilde{B}^\eps_1 +\tilde{B}^\eps_2
 \ee
 where 
 \be
 \begin{split}
 \tilde{B}^\eps_1=&  \int_{\partial B(0,1)} \left(\partial_\nu \tue \right)^2 +\frac{1}{2} \tue \partial_\nu \tue  - \frac{\vert \nabla \tue\vert^2}{2}  d\sigma \hbox{ and }\\
\tilde{B}^\eps_2 =& \int_{\partial B(0,1)}  \left(\frac{\mu_\eps}{r_\eps} \right)^2  \frac{\tue^6}{6}  d\sigma \hskip.1cm.
 \end{split}
 \ee 
Thanks to Claim \ref{estim2} and to Lebesgue dominated convergence theorem, we can pass to the limit to obtain that 
 \be
  \int_{\partial B(0,1)} \left(\partial_\nu U \right)^2 + \frac{1}{2} U\partial_\nu U - \frac{\vert \nabla U\vert^2}{2}  d\sigma =0 \hskip.1cm.
 \ee
Since $b$ is harmonic, it is easily checked that the left-hand side is just $-\frac{\omega_2 b(0)}{2}$. This proves that $b(0)=0$.

In order to prove the second part of the Claim, we apply the Poho\v{z}aev identity (\ref{pohozaev4}) of appendix \ref{Pohozaevidentities} to $\tue$ in $B(0,1)$. We obtain that 
\beq
\label{eq1.26}
\begin{split}
&\int_{\partial B(0,1)} \left( \frac{\vert \nabla \tue \vert^2}{2} \nu  -  \partial_\nu \tue \nabla \tue \right) \,d\sigma \\
&\quad = - \int_{B(0,1)} r_\eps^2 \tilde{h}_\eps \frac{\nabla \tue^2 }{2} \, dx-   \int_{\partial B(0,1)} \left(\frac{\mu_\eps}{r_\eps} \right)^2 \frac{ \tue^6}{6} \nu \, d\sigma \hskip.1cm.
\end{split}
\eeq 
It is clear that
\be
\int_{\partial B(0,1)} \left( \frac{\vert \nabla \tue \vert^2}{2} \nu -\partial_\nu \tue \nabla \tue  \right) d\sigma  \rightarrow \int_{\partial B(0,1)} \left(\frac{\vert \nabla U \vert^2}{2} \nu - \partial_\nu U \nabla U   \right) d\sigma \hbox{ as } \eps\rightarrow 0\hskip.1cm.
\ee
Moreover, thanks to the fact that $b$ is harmonic, we easily get that
\be
\int_{\partial B(0,1)} \left( \frac{\vert \nabla U \vert^2}{2} \nu - \nabla U \partial_\nu U\right)\, d\sigma = \omega_2 \nabla b(0)\hskip.1cm.
\ee
It remains to deals with the right-hand side of (\ref{eq1.26}). It is clear that
\be
\int_{\partial B(0,1)} \left(\frac{\mu_\eps}{r_\eps} \right)^2 \frac{\tue^6}{6}\nu d\sigma \rightarrow 0 \hbox{ as } \eps \rightarrow 0\hskip.1cm.
\ee 
Then we rewrite the first term of the right-hand side of (\ref{eq1.26}) as
\be
 \int_{B(0,1)} r_\eps^2 \tilde{h}_\eps \frac{\nabla \tue^2 }{2} \,dx =  \int_{ B(0,1)} r_\eps^2 (\tilde{h}_\eps - \tilde{h}_\eps(0)) \frac{\nabla \tue^2 }{2}\, dx + \tilde{h}_\eps(0) \int_{ B(0,1)} r_\eps^2  \frac{\nabla \tue^2 }{2} \,dx\hskip.1cm.
\ee
If we assume that the convergence of  $(h_\eps)$ holds in $C^{0,\eta}$, we can use Lebesgue dominated convergence theorem to obtain that the first term of the right-hand side goes to $0$ as $\eps\rightarrow 0$. Then, integrating by parts the second term, we get
\be
\tilde{h}_\eps(0) \int_{ B(0,1)} r_\eps^2  \frac{\nabla \tue^2 }{2} dx =  \tilde{h}_\eps(0) \int_{\partial B(0,1)} r_\eps^2  \frac{ \tue^2 }{2} \nu d\sigma
\ee
which clearly goes to $0$ as $\eps \rightarrow 0$. Finally, collecting the above informations, and passing to the limit $\eps\to 0$ in (\ref{eq1.26}), we get that $\nabla b(0)=0$ if the convergence of $(h_\eps)$ holds in $C^{0,\eta}$, which achieves the proof of the Claim.
\hfill$\blacksquare$

\medskip We are now in position to end the proof of proposition \ref{estim}.  If $\rho_\eps \rightarrow 0$ as $\eps\rightarrow 0$ then we deduce the proposition from claims \ref{estim3} and \ref{estim4}. If $\rho_\eps\not\to 0$ as $\eps\to 0$, then claims \ref{estim3} and \ref{estim4} give that $r_\eps\not \to 0$ as $\eps\to 0$. Then, using the Harnack inequality (\ref{eq1.16}), one can extend the result of Claim \ref{estim2} to $B(x_\eps,2\rho_\eps)\setminus\{x_\eps\}$, which proves the first part of Proposition \ref{estim} when $\rho_{\eps} \not\rightarrow 0$, and ends the proof of the whole proposition.

\section{Stability of the Poho\v{z}aev obstruction}\label{sectionstability}

We prove theorem \ref{thmstability} and give some stability result for radial solutions on the unit ball (see the end of the section). We assume by contradiction that there exists a sequence $\left(h_\eps\right)$ of functions  in $C^{0,\eta}\left(\R^3\right)$ for some $\eta>0$ and a sequence $\left(u_\eps\right)$ of $C^2$-solutions of (\ref{eq1.1}) where $\Omega$ is some smooth domain of $\R^3$ star-shaped with respect to the origin and $h_\eps \rightarrow h \hbox{ in } L^{p}\left(\Omega\right)$ as $\eps \rightarrow 0$ for some $p>3$ where $h\in C^1(\R^3)$ satisfies (\ref{mainhyp}). Sometimes  we will assume that $h_\eps \rightarrow h$ in $C^{0,\eta}$ as $\eps\rightarrow 0$. 

\medskip We claim first that 
\beq
\label{eq2.1}
\Vert \ue \Vert_\infty \rightarrow +\infty \hbox{ as }\eps \rightarrow 0\hskip.1cm.
\eeq
Indeed, if $(\ue)$ is uniformly bounded in $L^{\infty}\left(\Omega\right)$, then it is clear that $\left( \frac{\ue}{\Vert \ue \Vert_{\infty}}\right)$ is uniformly bounded in $W^{2,p}\left(\Omega\right)$ for some $p>3$, and thus, after passing to a subsequence, $\frac{\ue}{\Vert \ue \Vert_{\infty}} \rightarrow u $ in  $C^1_{loc}(\Omega)$ where $u$ is a positive solution of
\be
\Delta u + h u = \left( \lim_{\eps \rightarrow 0} \Vert \ue \Vert_{\infty}^{4} \right) u^5 \hbox{ in } \Omega 
\ee
with $u = 0$ on $\partial\Omega$. Since $h\geq 0$, it is clear that  $\Vert \ue \Vert_{\infty} \not\rightarrow 0$ as $\eps \rightarrow 0$. Then $\tilde{u}= \left( \lim_{\eps \rightarrow 0} \Vert \ue \Vert_{\infty} \right) u$ is a non-trivial solution of (\ref{mainequation}), which is a contradiction since (\ref{mainhyp}) holds. Thus (\ref{eq2.1}) is proved.

\medskip Then the sequence $(\ue)$ develops some concentration phenomena. We prove that this leads to a contradiction as follows~: in Claim \ref{iso1}, mimicking \cite{DruetHebey08}, we exhaust  a family of critical points of $\ue$, $\left(x_{1,\eps},\dots, x_{N_\eps,\eps}\right)$, such that each sequence $\left(x_{i_\eps, \eps}\right)$ satisfies the assumptions of Section \ref{sectionblowupanalysis} with
$$\rho_\eps= \min_{1\leq i\leq N_{\eps}, i\not=i_\eps} \left\{ \vert x_{i,\eps} -x_{i_\eps, \eps} \vert, d(x_{i_\eps, \eps}, \partial\Omega) \right\}\hskip.1cm.$$
In Claim \ref{iso4}, we prove that these concentration points are in fact isolated. In other words, we prove that $(\ue)$ develops only finitely many concentration points. We prove that such a configuration of concentrations points must satisfy two relations involving the Green function of $\Delta + h$ at these points. And it is impossible to find such a configuration thanks to some Poho\v{z}aev identity on Green functions we prove in Appendix \ref{pohog}. Claim \ref{iso1} is rather classical. The core of the proof lies in Claim \ref{iso4}. Avoiding bubble accumulation in the interior $\Omega$ in dimension $3$ is by now classical. The main difficulty here is to avoid boundary bubble accumulation. The rest of the section is devoted to the details of the proof we just sketched. 

\begin{claim}
\label{iso1} 
There exists $ D >0$ such that for all $\eps >0$, there exists $N_\eps \in {\mathbb N}^{*}$ and $N_\eps$ critical points of $\ue$, denoted by $(x_{1,\eps},\dots, x_{N_\eps})$ such that :
\be
\begin{split}
&d(x_{i,\eps},\partial\Omega ) \ue(x_{i,\eps})^2  \geq 1 \hbox{ for all } i\in [1,N_{\eps}]\hskip.1cm,\\
&\vert x_{i,\eps}-x_{j,\eps}\vert \ue(x_{i,\eps})^2  \geq 1 \hbox{ for all } i \not= j \in  [1,N_{\eps}]\hskip.1cm, 
\end{split}
\ee
and 
\be
\left( \min_{i\in [1,N_\eps ]} \vert x_{i,\eps} -x \vert  \right) \ue(x)^2 \leq D
\ee
for all $x\in \Omega$ and all $\eps>0$.
\end{claim}

\medskip {\bf Proof of Claim \ref{iso1}.} First of all, we claim that 
\beq
\label{eq2.2}
\{ x \in\Omega \hbox{ s.t. } \nabla \ue(x)=0 \hbox{ and } d(x,\partial\Omega) \ue(x)^2\geq 1 \}\not=\emptyset
\eeq
for $\eps$ small enough. Let us prove (\ref{eq2.2}). Let $y_\eps \in \Omega$ be a point where $\ue$ achieves his maximum. We set $\mu_\eps= \ue(y_\eps)^{-2}\rightarrow 0$ as $\eps\rightarrow 0$. We set also  for all $x\in\Omega_\eps= \{x\in\R^3 \hbox{ s.t. } y_\eps +\mu_\eps x \in\Omega \}$,
\be
\tue(x) =\mu_\eps^\frac{1}{2} \ue(y_\eps +\mu_\eps x)
\ee
which verifies
 \be
 \Delta \tue +\mu_\eps^2 \tilde{h}_\eps \tue = \tue^5 \hbox{ in } \Omega_{\eps},
 \ee
 where $\tilde{h}_\eps=h(y_\eps + \mu_\eps x)$. Note that $0\leq \tue\leq \tue(0)= 1$. Thanks to standard elliptic theory, we get that $\tue \rightarrow U$ in $C^{1}_{loc}(\Omega_0)$ where $U$ satisfies
\be 
\Delta U = U^5 \hbox{ in } \Omega_0 \hbox{ and } 0\leq U \leq 1=U(0),
\ee
and where $\Omega_0 =\lim_{\eps \rightarrow 0} \Omega_\eps $. Thanks to \cite{Dancer}, we have $\Omega_0=\R^3$, which proves that $d(y_\eps,\partial\Omega)\ue(y_\eps)^2 \rightarrow +\infty $ as $\eps\to 0$. This ends the proof of (\ref{eq2.2}).

Now, applying Lemma \ref{iso0}, see Appendix \ref{generallemma}, for $\eps$ small enough, there exist $N_\eps \in \N^{*}$ and $N_\eps$ critical points of $\ue$, denoted by $(x_{1,\eps},\dots, x_{N_\eps})$, such that :
\be
\begin{split}
&d(x_{i,\eps},\partial\Omega ) \ue(x_{i,\eps})^2  \geq 1 \hbox{ for all } i\in [1,N_{\eps}]\hskip.1cm,\\
&\vert x_{i,\eps}-x_{j,\vert} \ue(x_{i,\eps})^2  \geq 1 \hbox{ for all } i \not= j \in  [1,N_{\eps}] \hskip.1cm,
\end{split}
\ee
and 
\beq
\label{eq2.3}
\left( \min_{i\in [1,N_\eps ]} \vert x_{i,\eps} -x \vert  \right) \ue(x)^2 \leq 1
\eeq
for all critical point $x$ of $\ue$ such that $d(x,\partial\Omega)\ue(x)^2 \geq 1$. It remains to show that there exists $D>0$ such that 
\be
\left( \min_{i\in [1,N_\eps ]} \vert x_{i,\eps} -x \vert  \right) \ue(x)^2 \leq D
\ee
for all $x\in \Omega$. We proceed by contradiction, assuming that 
\beq
\label{eq2.4}
\sup_{x\in\Omega} \left( \left( \min_{i\in [1,N_\eps ]} \vert x_{i,\eps} -x \vert  \right) \ue^2(x)\right) \rightarrow +\infty
\eeq
as $\eps \rightarrow 0$. Let $z_\eps\in \Omega$ be such that
\be
\left( \min_{i\in [1,N_\eps ]} \vert x_{i,\eps} -z_\eps \vert  \right) \ue(z_\eps)^2 = 
\sup_{x\in\Omega} \left( \left( \min_{i\in [1,N_\eps ]} \vert x_{i,\eps} -x \vert  \right) \ue(x)^2\right) .
\ee
We set $\hat{\mu}_\eps =\ue(z_\eps)^{-2}$ and $S_\eps =\{ x_{1,\eps}, \dots, x_{N_\eps,\eps} \}$. Thanks to (\ref{eq2.4}), we check that 
\be 
\hat{\mu}_\eps \rightarrow 0\hbox{ as }\eps\rightarrow 0
\ee
and that 
\beq
\label{eq2.5}
\frac{d(S_\eps, z_\eps)}{\hat{\mu}_\eps} \rightarrow +\infty \hbox{ as } \eps\rightarrow 0 \hskip.1cm.
\eeq
Then we set, for all $x\in \hat{\Omega}_{\eps}=\{ x\in \R^3 \hbox{ s.t. } z_\eps + \hat{\mu}_\eps x \in \Omega \}$, 
\be 
\hue(x) = \hat{\mu}_\eps^\frac{1}{2} \hue(z_\eps+\hat{\mu}_\eps x) 
\ee
which verifies
\be
\Delta \hue +\hat{\mu}_\eps^2 \hat{h}_\eps \hue = \hue^5 \hbox{ in } \Omega_{\eps}
\ee
where $\hat{h}_\eps=h(z_\eps + \hat{\mu}_\eps x)$. Note that $\hue(0)=1$ and also that 
\be
\lim_{\eps\rightarrow 0} \sup_{B(0,R)\cap \Omega_\eps} \hue = 1
\ee
for all $R>0$ thanks to (\ref{eq2.4}) and (\ref{eq2.5}). Standard elliptic theory gives then that $\hue \rightarrow \hat{U}$ in $C^{1}_{loc}(\hat{\Omega}_0)$ where $U$ satisfies
\be 
\Delta \hat{U} = \hat{U}^5 \hbox{ in } \hat{\Omega}_0 \hbox{ and } 0\leq \hat{U} \leq1=\hat{U}(0)
\ee
with $\ds \hat{\Omega}_0=\lim_{\eps \rightarrow 0} \hat{\Omega}_\eps $. As above, we deduce that $\hat{\Omega}_0=\R^3$, which gives that 
\beq
\label{eq2.6}
\lim_{\eps\rightarrow 0} d(z_\eps, \partial\Omega) \ue^2(z_\eps)\rightarrow +\infty \hskip.1cm.
\eeq
Moreover, thanks to \cite{CGS}, we know that
\be
\hat{U}(x) = \frac{1} { \left(1+ \frac{\vert x \vert^{2} }{3} \right)^{\frac{1}{2}} } \hskip.1cm.
\ee
Since $\hat{U}$ has a strict local maximum at $0$, there exists $\hat{x}_\eps$, a critical point of $\ue$, such that $\vert z_{\eps} - \hat{x}_\eps \vert  = o(\hat{\mu}_\eps)$ and $\hat{\mu}_\eps \ue (\hat{x}_\eps)^{2}\rightarrow 1$ as $\eps \rightarrow 0$. Thanks to (\ref{eq2.5}) and (\ref{eq2.6}), this contradicts (\ref{eq2.3}) and proves the Claim.\hfill$\blacksquare$

\medskip We define
\be
d_\eps = \min \left\{ d(x_{i, \eps},x_{j, \eps}), d(x_{i, \eps},\partial\Omega) \hbox{ s.t. } 1\leq i < j \leq N_\eps \right\}
\ee
and prove~: 

\begin{claim}
\label{iso4}
If the convergence of $h_\eps$ to $h$ holds in $C^{0,\eta}$, then there exists $d>0$ such that  $\de \geq d$.
\end{claim}

\medskip {\bf Proof of Claim \ref{iso4}.} Assume that $\de \rightarrow 0 $ as $\eps \rightarrow 0 $. There are two cases to consider~: either the distance between two critical points goes to $0$, or one of them goes to the boundary. In the first case, the arguments which lead to a contradiction follow closely  \cite{DruetIMRN}, but in the second case we have to be more precise looking at the "artificial" singularities created by the boundary.

Up to reordering the concentration points, we can assume that 
\be
\de=d(x_{1,\eps}, x_{2,\eps}) \hbox{ or } d(x_{1,\eps}, \partial\Omega) \hskip.1cm. 
\ee
For $x\in \Omega_\eps =\{x\in\R^3 \hbox{ s.t. } x_{1,\eps} + d_{\eps} x \in \Omega \}$, we set
\be
\tue (x)= \de^{\frac{1}{2}}  u_\eps(x_{1,\eps} + \de x) 
\ee
which verifies
\be
\Delta \tue +\de^2 \tilde{h}_\eps \tue =\tue ^5 \hbox{ in }  \Omega_\eps, 
\ee
 where $\tilde{h}_\eps= h(x_{1,\eps} + \de x)$. We have, up to a harmless rotation, 
\be
\lim_{\eps \rightarrow 0} \Omega_\eps =\Omega_0=   \R^3 \hbox{ or } ]-\infty;d[ \times \R^2 \hbox{ where } d\geq 1\hskip.1cm.
\ee
We also set
\be
\tilde{x}_{i,\eps} =\frac{x_{i,\eps} - x_{1,\eps}}{\de}\hskip.1cm.
\ee
We claim that, for any sequence $i_\eps \in [1, N_\eps]$ such that
\beq
\label{eq2.7} 
\tue(\tilde{x}_{i_\eps,\eps})=O(1) \hskip.1cm,
\eeq
we have that
\beq 
\label{eq2.8}
\sup_{B(\tilde{x}_{i_\eps,\eps},\frac{1}{2})} \tue =O(1) \hskip.1cm.
\eeq
Indeed, let $y_{\eps} \in\overline{B(\tilde{x}_{i_\eps,\eps},\frac{1}{2})}$ be such that $\displaystyle \sup_{B(\tilde{x}_{i_\eps,\eps},\frac{1}{2})} \tue = \tue( y_{\eps} ) $ and assume by contradiction that 
\beq
\label{2.9}
\tue( y_{\eps} )^2 \rightarrow +\infty \hbox{ as }\eps \rightarrow 0\hskip.1cm.
\eeq
Thanks to the definitions of $\de$, $y_\eps$ and the last assertion of Claim \ref{iso1}, we can  write that
\be
\vert \de (y_\eps-\tilde{x}_{i_\eps,\eps}) \vert \ue(x_{1,\eps} +\de y_\eps)^2 \leq D
\ee
so that 
\beq
\label{2.10}
\vert y_\eps - \tilde{x}_{i_\eps,\eps}\vert = o(1)\hskip.1cm.
\eeq
For $x \in B(0, \frac{1}{3 \hat{\mu}_\eps})$ and $\eps$ small enough, we set 
\be
\hue (x)= \hat{\mu}_\eps^{\frac{1}{2}}  \tue (y_\eps + \hat{\mu}_\eps x)
\ee
where $\hat{\mu}_\eps= \ue(y_\eps)^{-2}$. It satisfies
\be
\Delta \hue +(\hat{\mu}_\eps d_\eps)^2 \hat{h}_\eps \hue =\hue ^5 \hbox{ in }  B(0, \frac{1}{3 \hat{\mu}_\eps})\hbox{ and }\hue(0) = \sup_{B(0, \frac{1}{3 \hat{\mu}_\eps})} \hue =1
\ee
where $\hat{h}_\eps= \tilde{h}_\eps(y_\eps + \hat{\mu}_\eps x)$. Thanks to (\ref{2.9}), $B(0, \frac{1}{3 \hat{\mu}_\eps}) \rightarrow \R^3$ as $\eps \rightarrow +\infty$. Then $\left(\hue\right)$ is uniformly locally bounded and, by standard elliptic theory, $\hue$ converges to $\hat{U}$ in $C^{1}_{loc}(\R^3)$ where $\hat{U}$ satisfies
\be 
\Delta \hat{U} = \hat{U}^5 \hbox{ in } \R^3 \hbox{ and } 0\leq \hat{U} \leq1=\hat{U}(0) \hskip.1cm.
\ee
Thanks to \cite{CGS} and to the fact that  $\frac{\tilde{x}_{i_\eps,\eps} - y_{\eps}}{\hat{\mu}_\eps}$ is bounded, we can write that
\be
\liminf_{\eps\to 0} \frac{\tue (x_{i_\eps,\eps})}{\tue(y_\eps)} >0
\ee
which is a contradiction with (\ref{eq2.7}) and (\ref{2.9}), and achieves the proof of (\ref{eq2.8}).

\medskip For $R>0$, we set $S_{R,\eps} = \left\{ \tilde{x}_{i,\eps} \vert \tilde{x}_{i,\eps}\in B(0, R)\right\}$. Thanks to definition of $\de$, up to a subsequence, $S_{R,\eps} \rightarrow S_{R}$ as $\eps \rightarrow 0$, where $S_{R}$ is a not empty finite set, then up to performing a diagonal extraction, we can define the countable set 
\be 
S=\bigcup_{R>0} S_{R}\hskip.1cm.
\ee 
Thanks to the previous definition, we are ready to prove the following assertion~:
\beq
\label{eq2.11}
\forall\,  i_{\eps} \in [1,N_\eps ]\hbox{ s.t. } d(x_{i_\eps ,\eps}, x_{1,\eps}) = O(\de)\, ,\,\, \tue(\tilde{x}_{i_\eps ,\eps}) \rightarrow +\infty\hbox{ as }\eps\to 0\hskip.1cm.
\eeq
Assume that there exists $i_\eps$ such that $d(x_{i_\eps ,\eps}, x_{1,\eps}) = O(\de)$ with $\tue(\tilde{x}_{i_\eps ,\eps})$ bounded, then for all sequences $j_\eps$ such that $d(x_{j_\eps ,\eps}, x_{1,\eps}) = O(\de)$,  $\tue(\tilde{x}_{j_\eps ,\eps})$ is bounded. Indeed, if  there exists a sequence $j_\eps$ such that $d(x_{j_\eps ,\eps}, x_{1,\eps}) = O(\de)$  and $\tue(\tilde{x}_{j_\eps ,\eps})\rightarrow +\infty $ as $\eps\rightarrow 0$, thanks to  Claim \ref{iso1}, we can apply Proposition \ref{estim} with $x_\eps = \tilde{x}_{j_\eps,\eps}$ and $\rho_{\eps} = \frac{\de}{3}$. We obtain that up to a subsequence  $\tue\rightarrow 0$ in $C^{1}_{loc}(B(\tilde{x},\frac{2}{3}))\setminus \{ \tilde{x} \}$, where $\ds \tilde{x}= \lim_{\eps \rightarrow 0} \tilde{x}_{j_\eps,\eps }$. But $\left(\tue\right)$ is uniformly bounded in $B(\tilde{y},\frac{1}{2})$, where  $\ds \tilde{y}= \lim_{\eps \rightarrow 0} \tilde{x}_{i_\eps,\eps }$. We thus obtain thanks to Harnack's inequality that $\tue(\tilde{x}_{i_\eps ,\eps})\rightarrow 0 $ as $\eps\rightarrow 0$,  which is a contradiction with the first or the second  assertion of Claim \ref{iso1}.

Thus we have proved that for all sequence $j_\eps$ such that $d(x_{j_\eps ,\eps}, x_{1,\eps}) = O(\de)$, $\tue(\tilde{x}_{j_\eps ,\eps})$ is bounded, which proves that $\left(\tue\right)$ is uniformly bounded in a neighborhood of any finite subset of $S$. But thanks to Claim \ref{iso1}, $\tue$ is bounded in any compact subset of $\Omega_{0}\setminus S$. This clearly proves that $\tue$ is uniformly bounded on any compact of $\Omega_0$. Then, by standard elliptic theory, $\tue \rightarrow U$ in $C^{1}_{loc}(\Omega_0)$ as $\eps \rightarrow 0$, where $U$ is a nonnegative solution of
\be
\Delta U = U^5 \hbox{ in } \Omega_0 \hskip.1cm.
\ee 
But, thanks to the first or second assertion of Claim \ref{iso1}, we know that $U(0) \geq 1 $, hence we have necessarily that  $\Omega_0=\R^3$, and thus  $U$ possesses at least two critical points, namely $0$ and $\ds \check{x}_2 =\lim_{\eps \rightarrow 0}  \check{x}_{2,\eps}$. Thanks to \cite{CGS}, this is impossible. This ends the proof of (\ref{eq2.11}).

\medskip We are now going to consider two cases, depending on $\Omega_0$. 

\medskip {\it Case 1 : $\Omega_0= \R^3$ -} In this case, up to a subsequence, $d_\eps=d(x_{1,\eps},x_{2,\eps})$ and $\ds S=\{0,\tilde{x}_2=\lim_{\eps \rightarrow 0} \tilde{x}_{2,\eps}, \dots  \}$ contains at least two points. Applying Proposition \ref{estim} with $x_\eps = \tilde{x}_{i,\eps}$ and $\rho_{\eps} = \frac{\de}{3}$, we obtain that 
\be
\tue(0)\tue(x) \rightarrow H= \frac{1}{\vert x \vert} +\frac{\lambda_{2}}{\vert x-\tilde{x}_2  \vert}+ \tilde{b} \hbox{ in }  C^{1}_{loc} (\R^3 \setminus S) \hbox{ as } \eps\rightarrow 0
\ee
where $\tilde{b}$ is an harmonic function in $\Omega_{0}\setminus \{S\setminus \{0, \check{x}_2 \}\}$, and $\lambda_2 >0$. Moreover $\tilde{b}(0)=-\lambda_2$. We prove in the following that $\tilde{b}$ is nonnegative, which will give a contradiction and end the study of this case. To check that $\tilde{b}$ is nonnegative, for all positive number $r$, we rewrite $H$ as 
\be
H= \sum_{\tilde{x}_i \in S\cap B(0,r)} \frac{\lambda_{i}}{\vert x-\tilde{x}_i \vert} + \hat{b}_r,
\ee
where $\lambda_i>0$. Then, taking $R>r$ big enough, we get  that  $\hat{b}_r > \frac{-1}{r}$ on $\partial B(0,R)$. Moreover, for any $\tilde{x}_j \in B(0,R)\setminus B(0,r)$, there exist a neighborhood $V_{j,r}$ of $\tilde{x}_j$ such that $ \hat{b}_r >0$ on $V_{j,r} $. Thanks to the maximum principle, $\hat{b}_r > \frac{-1}{r}$ on $B(0,R)$. Since $\hat{b}_r \rightarrow \hat{b}$ on every compact set as $r \rightarrow +\infty$, we get that $\ds H=  \sum_{\tilde{x}_i \in S} \frac{\lambda_{i}}{\vert x-\tilde{x}_i \vert} + \hat{b}$ with $\hat{b}\geq 0$, which proves that $\tilde{b}\geq 0$. This is the contradiction we were looking for, and this ends the proof of the claim in this first case.

\medskip {\it Case 2 : $\Omega_0= ]-\infty, d[\times \R^2$ -} We still denote $S=\{0=\tilde{x}_1, \tilde{x}_2, \dots \}$ and we apply Proposition \ref{estim} with $ x_\eps= x_{i,\eps}$ and $\rho_\eps= \frac{d_\eps}{3}$ to get that
\be
\tue(0)\tue(x) \rightarrow H= \sum_{\tilde{x}_i \in S}  \frac{\lambda_i}{\vert x-\tilde{x}_i \vert} + \tilde{b} \hbox{ in } C^{1}_{loc} (\Omega_0 \setminus S)
\ee
 where $\lambda_i  >0$, and $\tilde{b}$ is some harmonic function in $\Omega_0$. We extend $H$ to $\R^3$ by setting
$$
\hat{H}(x) = \left\{
\begin{array}{ll}
H(x) &\hbox{ if } x_1 \leq d \\
-H(s(x)) &\hbox{ otherwise }
\end{array}\right.
$$
where $s$ is the symmetry with respect to $ \{ d\} \times \R^2 $.  We also extend $\tilde{b}$ by setting  
\be
\hat{H} =  \sum_{\tilde{x}_i \in S}  \left( \frac{\lambda_i}{\vert x-\tilde{x}_i \vert} -  \frac{\lambda_i}{\vert s(x)-\tilde{x}_i\vert} \right) +\hat{b}\hskip.1cm.
\ee
It is clear that $\hat{b}$ is harmonic on $\R^3$ and satisfies $\hat{b} \geq 0$ in $\Omega_0$ and $\hat{b}\le 0$ in $\R^3\setminus \Omega_0$. This can be proved as in Case 1. Let $\mathcal{G}_R$ the Green function of the laplacian on the ball centered in $0$ with radius $R$, we get thanks to the Green representation formula that 
\be
\hat{b} (x) = \int_{\partial B(0,R)} \partial_\nu \mathcal{G}_R (x,y) \hat{b}(y) d\sigma 
\ee 
which gives since 
$$\partial_\nu {\mathcal G}_R\left(x,y\right)= \frac{R^2-\vert x\vert^2}{\omega_2 R \left\vert x-y\right\vert^3}$$
on $\partial B(0,R)$ that 
\be
\partial_1 \hat{b} (0) = \frac{3}{\omega_2 R^{4}} \int_{\partial B(0,R)} y_1 \hat{b}(y) d\sigma \hskip.1cm.
\ee 
Now we decompose $\partial B(0,R)$ into three sets, namely 
\begincal
A&=&\{y\in  \partial B(0,R)\hbox{ s.t. }  y_1 \geq d \}\hskip.1cm,\\ 
B&=&\{ y\in  \partial B(0,R) \hbox{ s.t. } 0 \leq y_1\leq d \}\hskip.1cm,\\
C&=&\{ y\in \partial B(0,R) \hbox{ s.t. } y_1\leq 0 \}\hskip.1cm.
\fincal 
In $A$ and $B$, we have that $y_1 \hat{b}(y) \leq d \hat{b}(y)$, and in $C$, we have that $y_1 \hat{b}(y) \leq 0$. Since $\hat{b}\ge 0$ in $C$, we arrive to 
\be
\partial_1 \hat{b} (0) \leq \frac{3d}{\omega_2 R^{4}} \int_{A\cup B}  \hat{b}(y) d\sigma \leq  \frac{3d}{\omega_2 R^{4}} \int_{\partial B(0,R)} \hat{b}(y) d\sigma = \frac{3d \hat{b}(0)}{R^2}\hskip.1cm.
\ee 
Passing to the limit $R\rightarrow +\infty$ gives that $\partial_1 \hat{b}(0) \leq 0$. In order to obtain a contradiction,  we rewrite $H$ in a neighborhood of $0$ as
 \be
 H(x) = \frac{1}{\vert x \vert} + \check{b}(x)
 \ee 
 where  
 \be
 \check{b}(x)= \hat{b}(x) -\frac{1}{\vert s(x) \vert } + \sum_{\check{x}_i \in S \setminus \{0\}} \lambda_i \left( \frac{1}{\vert x-\check{x}_i \vert} -  \frac{1}{\vert s(x) - \check{x}_i \vert} \right) \hskip.1cm. 
 \ee
As is easily checked, $\partial_1  \check{b}(0) <0$, which is a contradiction with Proposition \ref{estim}. This ends the proof of  Claim \ref{iso4} in this second case.

\medskip We are now ready to prove theorem \ref{thmstability}. Thanks to Claim \ref{iso1}, there exist $D >0, N \in \mathbb{N}^{*}$ and $N$ local maxima of $\ue$, $x_{1,\eps}, \dots, x_{N}$, such that:
\be
\begin{split}
&d( x_{i,\eps}, \partial \Omega )u_\eps\left(x_{i,\eps}\right)^2 \geq 1 \hbox{ for all } i\in[1,N]\hskip.1cm ,\\
&\vert x_{i,\eps}-x_{j,\eps}\vert \ue(x_{i,\eps})^2 \geq 1 \hbox{ for all } i \not= j\in[1,N]
\end{split}
\ee
and
\be
\left( \min_{i\in [1,N]} \vert x_{i,\eps}-x \vert \right) \ue(x)^2  \leq D
\ee
for all $x \in \Omega$. We can assume that $\ue(x_{i,\eps}) \rightarrow +\infty$ as $\eps \rightarrow 0$. Indeed, otherwise we can remove $x_{i,\eps}$ from the family of concentration points, and up to changing $D$, the assertion remains true. Then, thanks to Harnack inequality, there exists $C>0$ such that
\beq
\label{eq2.12}
\frac{1}{C} \ue(x_{1,\eps}) \leq \ue(x_{i,\eps}) \leq C \ue(x_{1,\eps})\hskip.1cm.
\eeq 
Now, thanks to the results of section \ref{sectionblowupanalysis} and by standard elliptic theory, we have that, after passing to a subsequence,
\be
\ue(x_{1,\eps}) \ue(x) \rightarrow G \hbox{ in } C^2_{loc} (\Omega \setminus \{x_1, \dots, x_N \}) \hbox {as } \eps \rightarrow 0
\ee
where 
\be
G(x)= \sum_{i=1}^{N} \lambda_i \mathcal{G}_{h}(x_i,x)
\ee
with $\mathcal{G}_{h}$ the Green function of the limit operator $\Delta + h$ with Dirichlet boundary condition on $\Omega$. Thanks to (\ref{eq2.12}), we know that $\lambda_i >0$ for $1\le i\le N$. This can be rewritten as 
\beq
\label{eq2.13}
G(x)=\frac{\lambda_i}{\omega_2 \vert x -x_i \vert} + G_i (x)
\eeq
where $G_i$ is a continuous function on $\Omega \setminus \{x_1, \dots,x_{i-1}, x_{i+1}, \dots, x_N \}$. Thanks to lemma \ref{decomp}, we can write 
\beq\label{eq2.16}
G_i(x)=  G_i(x_i) + \frac{h(x_i)}{2\omega_2} \vert x-x_i \vert + \gamma_i(x)
\end{equation}
where $\gamma_i\in C^1(\Omega)$ and $\gamma_i(0)=0$. We claim that 
\beq\label{eq2.14}
G_i\left(x_i\right) =0 \hbox{ for all }1\le i\le N\hskip.1cm.
\eeq
In order to prove this, we apply the Poho\v{z}aev identity (\ref{pohozaev1}) to $\ue$ on the ball $B\left(x_{i,\eps},\delta\right)$ for some $\delta >0$ small enough. This gives 
\beq\label{eq2.15}
\begin{split}
&\frac{1}{2} \int_{B(x_{i,\eps},\delta)} \left( h_{\eps} \ue^2 +  h_{\eps} <x-x_{i,\eps},\nabla \ue^2> \right) dx \\
&\quad =  \int_{\partial B(x_{i,\eps},\delta)}  \left( \delta\left( \partial_\nu \ue\right)^2  - \delta \frac{\vert \nabla \ue\vert^2}{2}+\frac{1}{2} u_\eps\partial_\nu u_\eps +\frac{\delta}{6} \ue^6\right)  \, d\sigma \hskip.1cm.
\end{split}
\eeq 
Thanks to the fact that $h_\eps$ is bounded in $L^p(\R^3)$ for some $p>3$ and Proposition \ref{estim}., we get the uniform estimate
\be
\ue\left(x_{i,\eps}\right)^2\left\vert \frac{1}{2} \int_{B(x_{i,\eps},\delta)} \left( h_{\eps} \ue^2 +  h_{\eps} <x-x_{i,\eps},\nabla \ue^2> \right) dx \right\vert\leq e(\delta) 
\ee
where $e\in C^0(\R)$ with $e(0)=0$. Using (\ref{eq2.13}), we get that 
\be
\begin{split}
& \int_{\partial B(x_{i,\eps},\delta)}  \left( \delta\left( \partial_\nu \ue\right)^2  - \delta \frac{\vert \nabla \ue\vert^2}{2}+\frac{1}{2} u_\eps\partial_\nu u_\eps \right)  \, d\sigma + \int_{\partial B(x_{i,\eps},\delta)}  \frac{\delta}{6} \ue^6\, d\sigma\\
&\quad = u_\eps\left(x_{i,\eps}\right)^{-2}\int_{\partial B(x_i,\delta)} \left( \delta \left( \partial_\nu G \right)^2 -  \delta \frac{\vert \nabla G \vert^2}{2}+\frac{1}{2} G\partial_\nu G \right)\, d\sigma +o\left(u_\eps\left(x_{i,\eps}\right)^{-2}\right)\hskip.1cm.
\end{split}
\ee
Using (\ref{eq2.16}), we easily get that  
\be
\int_{\partial B(x_i,\delta)} \left( \delta \left( \partial_\nu G \right)^2 -  \delta \frac{\vert \nabla G \vert^2}{2}+\frac{1}{2} G\partial_\nu G \right)\, d\sigma = -\frac{1}{2}\lambda_i G_i\left(x_i\right) +o(1) \hbox{ as } \delta\rightarrow 0\hskip.1cm.
\ee
Collecting the above informations proves (\ref{eq2.14}). 

\medskip We are going to prove now that $\nabla \gamma_i\left(x_i\right) = 0$ where $\gamma_i$ is as in (\ref{eq2.16}). This will contradict lemma \ref{pohoglemma} of Appendix \ref{pohog} and will achieve the proof of the theorem. For that purpose, we apply the Poho\v{z}aev identity (\ref{pohozaev4}) to $\ue$ on the ball  
$B\left(x_{i,\eps},\delta\right)$ for some $\delta >0$ small enough. We obtain that 
\beq\label{eq2.17}
\begin{split}
&\ue\left(x_{i,\eps}\right)^2\int_{\partial B(x_{i,\eps} ,\delta)} \left( \frac{\vert \nabla \ue \vert}{2} \nu -  \nabla \ue \partial_\nu \ue\right) d\sigma \\
&\quad = \ue\left(x_{i,\eps}\right)^2\int_{B(x_{i,\eps} ,\delta)} h_\eps \frac{\nabla \ue^2 }{2} dx -\ue\left(x_{i,\eps}\right)^2 \int_{\partial B(x_{i,\eps} ,\delta)} \nabla \ue^6\, d\sigma \hskip.1cm.
\end{split}
\eeq 
It is clear that we can pass to the limit in the left-hand side. Moreover, thanks to (\ref{eq2.14}) and (\ref{eq2.16}), we have that
\be
\int_{\partial B(x_i ,\delta)} \left( \frac{\vert \nabla G \vert}{2} \nu - \nabla G \partial_\nu G\right) d\sigma \rightarrow \nabla\gamma_{i}(x_i) \hbox{ as } \delta \rightarrow 0 .
\ee
Now we look at  the right-hand side of (\ref{eq2.17}). It is clear that 
$$\ue\left(x_{i,\eps}\right)^2 \int_{\partial B(x_{i,\eps} ,\delta)} \nabla \ue^6\, d\sigma\to 0\hbox{ as }\eps\to 0\hskip.1cm.$$ 
Then we write that 
\be
 \int_{B(x_{i,\eps},\delta)}  h_\eps \frac{\nabla \ue^2 }{2} \,dx =  \int_{ B(x_{i,\eps} ,\delta)} (h_\eps - h_\eps(x_{i,\eps})) \frac{\nabla \ue^2 }{2}\,dx + h_\eps(x_{i,\eps}) \int_{ B(x_{i,\eps},\delta)} \frac{\nabla \ue^2 }{2}\, dx\hskip.1cm.
\ee
Assuming that the convergence of $h_\eps$ to $h$ holds in $C^{0,\eta}$, it is clear that the first term of the right-hand side goes to $0$ as $\eps\rightarrow 0$. Integrating by parts the second term, we get that 
\be
h_\eps(x_{i,\eps}) \int_{ B(x_{i,\eps},\delta)}  \frac{\nabla \ue^2 }{2} dx =  h_\eps(x_{i,\eps}) \int_{\partial B(x_{i,\eps},\delta)} \frac{ \ue^2 }{2} \nu \, d\sigma\rightarrow  h\left(x_i\right) \int_{\partial B(x_i,\delta)} \frac{G^2}{2} \nu d\sigma
\ee 
as $\eps\to 0$. It is easily checked that the above goes to $0$ as $\delta\rightarrow 0$.

Finally, collecting the above informations, and passing consecutively to the limit $\eps\to 0$ and $\delta\to 0$ in  (\ref{eq2.17}), we get that $\nabla\gamma_{i}(x_i)=0$ for all $i$, which achieves the proof of theorem \ref{thmstability} thanks to lemma \ref{pohog}.\hfill $\blacksquare$

\medskip Let us now give a precise statement of what we meant by stability of the Poho\v{z}aev obstruction in the radial situation in the introduction. We will prove the following~:

\begin{thm}\label{propradialstability}
Let $B$ be the unit ball of $\R^3$. Let $h_0$ be a $C^1$-radial function which satisfies (\ref{mainhyp}). Then for any $p>3$, there exists $\delta>0$ (depending on $h_0$ and $p$) such that if $h\in C^{0,\eta}\left(B\right)$ for some $\eta>0$ with $\left\Vert h-h_0\right\Vert_{L^p(B)} \le \delta$, then there exists no positive radial solution of equation (\ref{mainequation}) in the unit ball. 
\end{thm}

\medskip {\bf Proof of theorem \ref{propradialstability}.} We proceed as for the proof of theorem \ref{thmstability}. Note that, since $u_\eps$ is radial, there can be only one concentration point, namely $0$. Thanks to claim \ref{iso1},  the result of section \ref{sectionblowupanalysis} and standard elliptic theory, we have that 
\be
\ue(0) \ue(x) \rightarrow \omega_2 \mathcal{G}_{h}(x,0) \hbox{ in } C^1_{loc} (\Omega\setminus\{0\}) \hbox {as } \eps \rightarrow 0
\ee
where $\mathcal{G}_{h}$ is the Green function of the limit operator $\Delta + h$. We have that 
\be
\mathcal{G}_{h}(x,0)= \frac{1}{\omega_{2} \vert x \vert} +g(x)
\ee
where $g$ is a continuous function on $\Omega$ which satisfies 
\be
\Delta g + hg = -\frac{h}{\omega_2 \vert x \vert} \hbox{ in } \Omega\hbox{ and }g =-\omega_2 \hbox{ on } \partial\Omega\hskip.1cm.
\ee
By the maximum principle, we see that $g$ is negative so that $g(0)<0$. Now we can proceed as in the proof of (\ref{eq2.14}) to get a contradiction. Note that the proof of (\ref{eq2.14}) did not require the $C^{0,\eta}$ convergence of $h_\eps$. In the above proof, this $C^{0,\eta}$ convergence had been used only in the proof of claim \ref{iso4} (which is given for free in the radial situation) and in the 
last part of the proof to deal with the case of several concentration points (this can not happen in the radial situation). \hfill  $\blacksquare$

\medskip We shall prove in the next section that the above theorem is sharp in the radial situation.

\section{Construction of blowing-up examples and instability of the Poho\v{z}aev obstruction}\label{sectionexamples}

In this section, we prove theorem \ref{thminstability}. In fact, we will first prove the corresponding result in the radial situation (thus showing that our theorem \ref{propradialstability} is sharp) since it contains the main ideas and the computations are a little bit less involved.

We first need some results on Green's functions of coercive operators $\Delta + h$ with Dirichlet boundary condition on domains of ${\mathbb R}^3$. We let $\Om$ be a smooth domain of ${\mathbb R}^3$ and $h\in C^1\left(\Omega\right)$ be such that the operator $\Delta + h$, with Dirichlet boundary conditions, is coercive. Then there exists a unique function ${\mathcal G}: \Om\times\Om \setminus \left\{\left(x,x\right),\, x\in \Om\right\}\mapsto {\mathbb R}$, symmetric, positive, such that 
$$\Delta_y {\mathcal G}\left(x,y\right) + h(y) {\mathcal G}\left(x,y\right) = \omega_2 \delta_x$$
in $\Om$ and ${\mathcal G}\left(x,y\right) =0$ for $y\in \partial\Om$ for all $x\in \Om$. It is easily checked that ${\mathcal G}\left(x,y\right)$ has the following expansion in the neighbourhood of the diagonal 
\begin{equation}\label{eq3.1}
{\mathcal G}\left(x,y\right) =\frac{1}{\left\vert x-y\right\vert} + \frac{1}{2} h(x) \left\vert x-y\right\vert + \gamma_x(y)
\end{equation}
where $\gamma_x\in C^1\left(\Omega\right)$ satisfies 
$$\Delta_y \gamma_x(y) + h(y) \gamma_x(y)= \frac{h(x)-h(y)}{\left\vert x-y\right\vert}-\frac{1}{2} h(x)h(y) \left\vert x-y\right\vert$$
in $\Om$ with 
$$\gamma_x (y) = -\frac{1}{\left\vert x-y\right\vert} - \frac{1}{2} h(x) \left\vert x-y\right\vert$$
for all $y\in \partial\Om$.

\subsection{The radial case} We start by proving that the Poho{\v{z}}aev identity is not $L^3_r$-stable in the unit ball of ${\mathbb R}^3$. More precisely, we prove the following result~: 

\begin{thm}\label{propinstabilityradial}
Let $h\in C^1\left(B\right)$ be a non-negative radial function on the unit ball $B$ of ${\mathbb R}^3$. For any $\eps>0$, there exists a radial function $\tilde{h}\in C^{0,\eta}\left(B\right)$ with $\left\Vert \tilde{h}-h\right\Vert_{L^3\left(B\right)} \le \eps$ such that the equation 
$$\Delta \tilde{u} + \tilde{h} \tilde{u} = \tilde{u}^5\hbox{ in }B\, \,\, \tilde{u}=0\hbox{ on }\partial B$$
admits a positive radial solution.
\end{thm}

\medskip Note that a function $h$ which satisfies (\ref{mainhyp}) is necessarily non-negative. Let us prove the theorem in the rest of this subsection. We let $h\in C^1\left(B\right)$ be a non-negative radial function where $B$ is the unit ball of ${\mathbb R}^3$. We let ${\mathcal G}$ be the Green function of $\Delta + h$ and let $G(x)={\mathcal G}\left(0,x\right)$. We set 
\begin{equation}\label{eq3.2}
u_\eps(x) = U_\eps(x) + \eta_\eps(x) V_\eps(x)
\end{equation}
where 
\begin{equation}\label{eq3.3}
\begin{split}
&U_\eps(x) = \eps^{\frac{1}{2}}\left(\eps^2+G(x)^{-2}\right)^{-\frac{1}{2}}\,\\
&V_\eps(x) = - \gamma_0(0)\eps^{\frac{1}{2}} G(x)^{-3}\left(\eps^2 + G(x)^{-2}\right)^{-\frac{3}{2}}\,\\
&\eta_\eps(x) =\eta(x) \frac{\ln \left(\eps^2+\vert x\vert^2\right)}{\ln \eps^2}
\end{split}
\end{equation}
where $\eta$ is a smooth positive function such that $\eta=1$ on the ball of radius $\frac{1}{4}$ and $\eta=0$ outside of the ball of radius $\frac{1}{2}$. Here, $\gamma_0$ comes from the asymptotic expansion (\ref{eq3.1}). It is easily checked that $u_\eps$ is a $C^{2,\eta}$-positive function in $B$ and that $u_\eps=0$ on $\partial B$. Moreover, we have that 
\begin{equation}\label{eq3.4}
\frac{\eta_\eps V_\eps}{U_\eps}\to 0 \hbox{ in }L^\infty\left(B\right)\hbox{ as }\eps\to 0\hskip.1cm.
\end{equation}
We claim that 
\begin{equation}\label{eq3.5}
\frac{3 u_\eps^5 - \Delta u_\eps}{u_\eps}\to h\hbox{ in }L^3\left(B\right)\hbox{ as }\eps\to 0\hskip.1cm,
\end{equation}
which clearly implies the theorem. Straightforward computations give that 
\begin{equation}\label{eq3.6}
\Delta U_\eps + h U_\eps = 3 U_\eps^5 \left\vert \nabla G^{-1}\right\vert^2+ h(x) \eps U_\eps^3
\end{equation}
and that 
\begin{equation}\label{eq3.7}
\begin{split}
\Delta V_\eps + h V_\eps = \, &15 U_\eps^4 V_\eps +12\gamma_0(0) \eps^{\frac{5}{2}} G^4 \left(1+\eps^2G^2\right)^{-\frac{5}{2}}\\
& - \eps^{\frac{1}{2}} \gamma_0(0) h \left(1+\eps^2G^2\right)^{-\frac{5}{2}} \left(1+4\eps^2 G^2\right)\\
& - 3\gamma_0(0) \eps^{\frac{5}{2}} G^4 \left(1+\eps^2 G^2\right)^{-\frac{7}{2}} \left(1-4\eps^2 G^2\right)\left(\left\vert \nabla G^{-1}\right\vert^2-1\right)
\hskip.1cm.
\end{split}
\end{equation}
It is easily checked that this implies that 
\begin{equation}\label{eq3.8}
\Delta u_\eps + h u_\eps - 3u_\eps^5 = o\left(u_\eps\right)
\end{equation}
in $B_0(1)\setminus B_0\left(\frac{1}{2}\right)$. Using the expansion of $G$ and its consequence 
$$\left\vert \nabla G^{-1}\right\vert^2 = 1 - 4\gamma_0(0)G^{-1} + O\left(G^{-2}\right)\hskip.1cm,$$
we can then write thanks to (\ref{eq3.4}) that 
\begin{equation}\label{eq3.9}
\begin{split}
\frac{\Delta u_\eps + h u_\eps - 3 u_\eps^5}{u_\eps} =&\, O\left(\left\vert x\right\vert^2 U_\eps^4\right)+O\left(\eps U_\eps^2\right)+ O\left(\left\vert x\right\vert U_\eps^4 \left\vert 1- \eta_\eps\right\vert\right) \\
&\, +O\left(U_\eps^{-1}\left\vert \nabla V_\eps\right\vert \left\vert \nabla \eta_\eps\right\vert\right) + O\left(U_\eps^{-1}\left\vert V_\eps\right\vert \left\vert \Delta \eta_\eps\right\vert\right)
\end{split}
\end{equation}
in $B_0\left(\frac{1}{2}\right)$. It is easily checked that 
\begin{equation}\label{eq3.10}
\left\vert x\right\vert^2 U_\eps^4\to 0 \hbox{ and }\eps U_\eps^2\to 0\hbox{ in }L^p\left(B\right)\hbox{ as }\eps\to 0
\end{equation}
for all $1\le p<+\infty$. Let us write now that 
\begincal
\int_{B} \left\vert x\right\vert^3 U_\eps^{12} \left\vert 1- \eta_\eps\right\vert^3\, dx &= &
O\left(\eps^6 \int_0^1 r^5 \left(\eps^2 + r^2\right)^{-6} \left\vert 1- \frac{\ln\left(\eps^2 +r^2\right)}{\ln \eps^2}\right\vert^3\, dr\right)\\
&= & O\left(\int_0^{\eps^{-1}} r^5 \left(1 + r^2\right)^{-6} \left\vert \frac{\ln \left(1+r^2\right)}{\ln \eps^2}\right\vert^3\, dr\right)\\
&=& O\left(\left\vert \ln \eps^2\right\vert^{-3}\right)=o(1) 
\fincal
thanks to the dominated convergence theorem. We can also write that 
\begincal
\int_B U_\eps^{-3}\left\vert \nabla V_\eps\right\vert^3 \left\vert \nabla \eta_\eps\right\vert^3 \, dx &=& 
O\left(\left\vert \ln \eps^2\right\vert^{-3}\int_0^1 r^{11} \left(\eps^2+r^2\right)^{-6}\, dr\right) \\
&&+ O\left(\int_{\frac{1}{4}}^{\frac{1}{2}} \left\vert \frac{\ln \left(\eps^2 +r^2\right)}{\ln \eps^2}\right\vert^3\, dr\right)\\
&=&O\left(\left\vert \ln \eps^2\right\vert^{-3}\right)=o(1) 
\fincal
and that 
\begincal 
\int_B U_\eps^{-3}\left\vert V_\eps\right\vert^3 \left\vert \Delta \eta_\eps\right\vert^3&=& O\left(\left\vert \ln \eps^2\right\vert^{-3}\int_0^1 r^5 \left(\eps^2 + r^2\right)^{-1}\, dr \right) + O\left(\left\vert \ln \eps^2\right\vert^{-3}\right)\\
&=& O\left(\left\vert \ln \eps^2\right\vert^{-3}\right)=o(1)\hskip.1cm.
\fincal
Coming back to (\ref{eq3.9}) with these last estimates, we get (\ref{eq3.5}). This ends the proof of  theorem \ref{propinstabilityradial}. \hfill $\blacksquare$

\subsection{The general case} Here we prove that the Poho{\v{z}}aev identity is never $L^\infty$-stable. In fact we will even prove a stronger result~:

\begin{thm}\label{propexemples}
Let $\Omega$ be a smooth domain of $\mathbb{R}^3$ and let $h\in C^1\left(\Omega\right)$ be such that the operator $\Delta + h$ is coercive. For any $\eps>0$, there exists $\tilde{h}\in C^{0,\eta}\left(\Om\right)$ with $\left\Vert \tilde{h}- h\right\Vert_\infty \le \eps$ such that the equation
$$\left\{\begin{array}{l}
{\displaystyle \Delta \tilde{u} + \tilde{h} \tilde{u} = \tilde{u}^5\hbox{ in }\Omega}\\
\,\\
{\displaystyle \tilde{u}=0\hbox{ on }\partial\Omega\,,\,\, \tilde{u}>0\hbox{ in }\Omega}
\end{array}\right.$$
admits a solution. 
\end{thm}

\medskip It is clear that this result implies theorem \ref{thminstability}. It is sufficient to remember that a function $h$ which satisfies (\ref{mainhyp}) is necessarily non-negative and that a non-negative $h$ leads to a coercive operator $\Delta+h$. The rest of this subsection is devoted to the proof of this theorem. 

\medskip We will construct a sequence of functions $u_\eps\in C^\infty\left(\Omega\right)$, positive in $\Omega$, null on the boundary of $\Omega$, such that 
\begin{equation}\label{eq3.11} 
\frac{\Delta u_\eps - 3 u_\eps^5}{u_\eps}\to h \hbox{ in }L^\infty\left(\Omega\right)\hbox{ as }\eps\to 0\hskip.1cm.
\end{equation}
This will clearly prove the theorem. 

\medskip We let ${\mathcal G}$ be the Green function of the operator $\Delta + h$ in $\Om$ with Dirichlet boundary conditions. Note first that $\gamma_x(x)\to -\infty$ as $x$ approaches $\partial \Om$. In particular, there exists a point $x_1\in \Om$ such that $\gamma_{x_1}\left(x_1\right)<0$. For $x\in \Om\setminus \left\{x_1\right\}$, we set 
$$\lambda(x) = \left(-\frac{\gamma_{x_1}\left(x_1\right)}{{\mathcal G}\left(x_1,x\right)}\right)^2$$
and 
$$F(x) = {\mathcal G}\left(x_1,x\right)^2 - \gamma_{x_1}\left(x_1\right) \gamma_{x}\left(x\right)\hskip.1cm.$$
Since $F(x)\to +\infty$ as $x\to x_1$ and $F(x)\to -\infty$ as $x$ approaches $\partial \Om$ and since $F$ is continous, there exists $x_2$ such that $F(x_2)=0$. We let then $\lambda = \lambda\left(x_2\right)$ and we have 
\begin{equation}\label{eq3.12}
\sqrt{\lambda} G_2\left(x_1\right) + \gamma_1\left(x_1\right) = G_1\left(x_2\right)+\sqrt{\lambda} \gamma_2\left(x_2\right)=0
\end{equation}
where 
\begin{equation}\label{eq3.13}
G_1(x) = {\mathcal G}\left(x_1,x\right)\, ,\,\,
G_2(x) = {\mathcal G}\left(x_2,x\right)\, ,\,\,
\gamma_1(x)=\gamma_{x_1}\left(x\right)\, ,\,\,
\gamma_2(x)=\gamma_{x_2}\left(x\right)\hskip.1cm.
\end{equation}
We let $\delta >0$ be such that $\delta \le 10 d\left(x_1,\partial \Om\right)$ and $\delta \le 10 d\left(x_2,\partial \Om\right)$. We fix $\eta\in C^\infty\left({\mathbb R}\right)$ such that $\eta(r)=1$ for $\left\vert r\right\vert \le \delta$ and $\eta(r)=0$ for $\left\vert r\right\vert \ge 2\delta$. We set in the following 
\begin{equation}\label{eq3.14}
\begin{split}
u_\eps = &\, \eps^{-\frac{1}{2}} U \left(\eps G_1\right) + \left(\lambda \eps\right)^{-\frac{1}{2}}U\left(\lambda \eps G_2\right) \\
&\, + \eta \left(\left\vert x-x_1\right\vert\right) \gamma_1\left(x_1\right) \eps^{\frac{1}{2}} V\left(\eps G_1\right)+ \eta \left(\left\vert x-x_2\right\vert\right) \gamma_2\left(x_2\right) \left(\lambda\eps\right)^{\frac{1}{2}} V\left(\lambda\eps G_2\right) \\
&\, - \eta \left(\left\vert x-x_1\right\vert\right) \psi_\eps\left(\eps G_1\right)\eps^{\frac{1}{2}}\left(x-x_1\right)^i\\
&\qquad . \left(\left(1+\eps^2G_1^2\right)^{-\frac{3}{2}}\partial_i\gamma_1\left(x_1\right)+\lambda ^{\frac{1}{2}} \partial_i G_2\left(x_1\right)\right) \\
&\, - \eta \left(\left\vert x-x_2\right\vert\right) \psi_\eps\left(\lambda\eps G_2\right)\left(\lambda \eps\right)^{\frac{1}{2}}\left(x-x_2\right)^i\\
&\qquad . \left(\left(1+\lambda^2\eps^2G_2^2\right)^{-\frac{3}{2}}\partial_i\gamma_2\left(x_2\right)+\lambda ^{-\frac{1}{2}} \partial_i G_1\left(x_2\right)\right) \\
&\, +  \eta \left(\left\vert x-x_1\right\vert\right) \eps^{\frac{3}{2}} \psi_\eps\left(\eps G_1\right)
\left(h\left(x_1\right)W\left(\eps G_1\right)-\frac{3}{2}\gamma_1\left(x_1\right)^2 U \left(\eps G_1\right)^5\right) \\
&\, +  \eta \left(\left\vert x-x_2\right\vert\right) \left(\lambda\eps\right)^{\frac{3}{2}} \psi_\eps\left(\lambda\eps G_2\right)
\left(h\left(x_2\right)W\left(\lambda\eps G_2\right)-\frac{3}{2}\gamma_2\left(x_2\right)^2 U \left(\lambda\eps G_2\right)^5\right) \\
\end{split}\end{equation}
where we adopt Einstein summation conventions and $U$, $V$, $W$ and $\psi_\eps$ are given by 
\begin{equation}\label{eq3.15}
\begin{split}
&U(r)= r\left(1+r^2\right)^{-\frac{1}{2}},\, 
V(r)= 1-\left(1+r^2\right)^{-\frac{3}{2}},\,
\psi_\eps(r) = 1+ \frac{\ln \left(1+r^{-2}\right)}{\ln \eps^2}\hbox{ and }\\
&W(r) = -\frac{13}{4} U + 8\left(2U^3-U\right)\ln U - 2\left(U^{-1}-8U+8U^3 \right) r\arctan\left(\frac{1}{r}\right)\hskip.1cm.
\end{split}
\end{equation}
 It is easily checked that $u_\eps$ is $C^{2,\eta}$ in $\Om$ and that $u_\eps=0$ on $\partial\Om$. We claim now that (\ref{eq3.11}) holds for this specific $u_\eps$ and that $u_\eps$ is positive in $\Omega$. We shall prove this claim in three steps. First, we can prove it rather easily in $\Om\setminus B_{x_1}\left(2\delta\right)\bigcup B_{x_2}\left(2\delta\right)$ because, in this region, $u_\eps$ is simply 
$$u_\eps = \eps^{-\frac{1}{2}} U \left(\eps G_1\right) + \left(\lambda \eps\right)^{-\frac{1}{2}}U\left(\lambda \eps G_2\right) \hskip.1cm.$$
Now, noticing that $U'=r^{-3}U^3$ and that $U'' = - 3 r^{-4}U^5$, simple computations lead to 
\begin{equation}\label{eq3.20}
\begin{split}
&\Delta \left(\eps^{-\frac{1}{2}} U \left(\eps G_1\right) + \left(\lambda \eps\right)^{-\frac{1}{2}}U\left(\lambda \eps G_2\right)\right) + h  \left(\eps^{-\frac{1}{2}} U \left(\eps G_1\right) + \left(\lambda \eps\right)^{-\frac{1}{2}}U\left(\lambda \eps G_2\right)\right)\\
&\quad = 3\left(\eps^{-\frac{1}{2}} U \left(\eps G_1\right)\right)^5 \left\vert \nabla G_1^{-1}\right\vert^2 
+3\left(\left(\lambda\eps\right)^{-\frac{1}{2}} U \left(\lambda\eps G_2\right)\right)^5 \left\vert \nabla G_2^{-1}\right\vert^2 \\
&\qquad + h \eps^{-\frac{1}{2}} U \left(\eps G_1\right) \left(1 - \left(1+\eps^2 G_1^2\right)^{-1}\right) \\
&\qquad +h  \left(\lambda \eps\right)^{-\frac{1}{2}}U\left(\lambda \eps G_2\right)\left(1 - \left(1+\lambda^2 \eps^2 G_2^2\right)^{-1}\right) 
\end{split}
\end{equation}
in $\Om$. In the region we are interested in, this clearly leads to 
$$\Delta u_\eps+h u_\eps-3u_\eps^5 = o\left(u_\eps\right)$$
which proves that (\ref{eq3.11}) holds in this region while $u_\eps$ is clearly positive in this region.

We will now prove that (\ref{eq3.11}) holds in $B_{x_1}\left(2\delta\right)$ and that $u_\eps$ is positive in this ball. By symmetry, it is clear that the proof of the fact that (\ref{eq3.11}) holds in $B_{x_2}\left(2\delta\right)$ is exactly the same\footnote{The symmetry is precisely the following~: if we see $u_\eps$ as a function of $x_1$, $x_2$, $\eps$ and $\lambda$, namely $u_\eps\left(x_1,x_2,\eps,\lambda\right)$, one has that $u_\eps\left(x_2,x_1,\lambda \eps,\lambda^{-1}\right)=u_\eps\left(x_1,x_2,\eps,\lambda\right)$.}. In order to simplify the notations, we will assume that $x_1=0$, which we can always do by translating $\Omega$. We will denote $G_1$ by $G$ and $\gamma_1$ by $\gamma$. We also set 
\begin{equation}\label{eq3.21}
\begin{split}
&U_\eps = \eps^{-\frac{1}{2}} U\left(\eps G\right),\, V_\eps = \eps^{\frac{1}{2}} V\left(\eps G\right),\, W_\eps = \eps^{\frac{3}{2}} W\left(\eps G\right),\\
&\varphi_\eps = \psi_\eps\left(\eps G\right),\, Y_\eps = -\frac{3}{2}\eps^{\frac{3}{2}} U\left(\eps G\right)^5,\, \tilde{U}_\eps = \left(\lambda\eps\right)^{-\frac{1}{2}} U\left(\lambda \eps G_2\right)\hbox{ and }\\
&Z_\eps = -\eps^{\frac{1}{2}}x^i\left(\left(1+\eps^2 G^2\right)^{-\frac{3}{2}}\partial_i\gamma\left(0\right)+\lambda^{\frac{1}{2}}\partial_i G_2\left(0\right)\right)\hskip.1cm.
\end{split}
\end{equation}
With these notations, we have that, in $B_0\left(2\delta\right)$, 
\begin{equation}\label{eq3.22}
u_\eps =  U_\eps + \tilde{U}_\eps +\eta\left(\vert x\vert\right) \gamma(0)V_\eps+\eta\left(\vert x\vert\right) \varphi_\eps \left(h(0)W_\eps + \gamma(0)^2 Y_\eps +Z_\eps\right)\hskip.1cm.
\end{equation}
Let us write thanks to (\ref{eq3.1}) that 
\begin{equation}\label{eq3.23}
\begin{split}
\left\vert \nabla G^{-1}\right\vert^2 = \, &\, 1-4\gamma(0)G^{-1} +3\left(2\gamma(0)^2-h(0)\right)G^{-2}\\
&\, - 6 G^{-1} x^i\partial_i\gamma(0)+o\left(G^{-2}\right)
\end{split}
\end{equation}
and let us also remark that 
\begin{equation}\label{eq3.23bis}
G^{-2} = \eps U_\eps^{-2} -\eps^2\hskip.1cm. 
\end{equation}
Let us write thanks to (\ref{eq3.12}) that 
\begin{equation}\label{eq3.24}
\begin{split}
&\tilde{U}_\eps =  -\eps^{\frac{1}{2}}\gamma(0)+\left(\lambda\eps\right)^{\frac{1}{2}} x^i \partial_i G_2(0) + O\left(\eps^{\frac{3}{2}} U_\eps^{-2}\right),\\ 
&V_\eps = O\left(\eps^{\frac{3}{2}}U_\eps^2\right), \, 
W_\eps =O\left(\eps^{\frac{3}{2}}\right),\, 
Y_\eps =O\left(\eps^{\frac{3}{2}}\right)\hbox{ and }Z_\eps =O\left(\eps U_\eps^{-1}\right)\hskip.1cm.
\end{split}
\end{equation}
Thanks to (\ref{eq3.22}) and to (\ref{eq3.24}), it is easily checked that $u_\eps$ is positive in $B_0\left(2\delta\right)$. Lengthy but straightforward computations lead then to
\begin{equation}\label{eq3.25}
\begin{split}
\left\vert \nabla \varphi_\eps\right\vert =\, &\, O\left(\frac{1}{\eps^{\frac{1}{2}}\ln \frac{1}{\eps}} U_\eps\right)\hskip.1cm,\hskip.1cm \left\vert \nabla V_\eps\right\vert =O\left(\eps U_\eps^3\right)\hskip.1cm,\hskip.1cm
\left\vert \nabla W_\eps\right\vert =  O\left(\eps U_\eps\right)\hskip.1cm,\\
\left\vert \nabla Y_\eps\right\vert =\, &\, O\left(\eps U_\eps\right)\hbox{ and } 
\left\vert \nabla Z_\eps\right\vert = O\left(\eps^{\frac{1}{2}}\right)\\
\end{split}
\end{equation}
in $B_0\left(2\delta\right)$ and to
\begin{equation}\label{eq3.26}
\begin{split}
\Delta U_\eps +h U_\eps = \, &\, 3U_\eps^5 - 12\gamma(0)G^{-1}U_\eps^5 - 18G^{-1} x^i\partial_i\gamma(0)U_\eps^5 \\
&\, + 18 \gamma(0)^2 G^{-2} U_\eps^5 + h(0)\left(\eps^2 - 8 G^{-2}\right) U_\eps^5 + o\left(U_\eps\right)\hskip.1cm,\\
\Delta \tilde{U}_\eps + h \tilde{U}_\eps = \,&\,  3\tilde{U}_\eps^5 \left\vert \nabla G_2^{-1}\right\vert^2 +\lambda\eps h \tilde{U}_\eps^3 = O\left(\eps^{\frac{5}{2}}\right)\hskip.1cm,\\
\Delta V_\eps + h V_\eps= \,&\, 15 U_\eps^4 V_\eps  - 15 \eps^{\frac{1}{2}} U_\eps^4+ 12 G^{-1}U_\eps^5 \\
&\, + 12 \gamma(0) \left(5\eps^{-1} G^{-4} U_\eps^7 - 4 G^{-2} U_\eps^5\right)+o\left(U_\eps\right)\hskip.1cm,\\
\Delta W_\eps +h W_\eps = \, &\, 15 U_\eps^4 W_\eps + 8\eps U_\eps^3 -9\eps^2 U_\eps^5 + o\left(U_\eps\right)\hskip.1cm,\\
\Delta Y_\eps +h Y_\eps = \, &\, 15 U_\eps^4 Y_\eps  +30\eps^3 U_\eps^7 -30 \eps^4 U_\eps^9+ o\left(U_\eps\right)\hskip.1cm,\\
\Delta \varphi_\eps =\, &\,O\left(\frac{1}{\eps \ln\frac{1}{\eps}} U_\eps^2\right)\hskip.1cm, \\
\Delta Z_\eps + h Z_\eps =\, &\, 15 U_\eps^4 Z_\eps + 18 U_\eps^5 G^{-1} \partial_i\gamma(0)x^i\\
&\, + 15 \left(\lambda \eps\right)^{\frac{1}{2}} U_\eps^4 x^i \partial_i G_2(0) + O\left(\eps U_\eps^{-1}\right)+o\left(U_\eps\right) \\
\end{split}
\end{equation}
in $B_0\left(2\delta\right)$. It follows easily from the above equations that 
$$\frac{\Delta u_\eps +h u_\eps - 3u_\eps^5}{u_\eps} \to 0 \hbox{ in }L^\infty\left(B_0\left(2\delta\right)\setminus B_0\left(\delta\right)\right)\hbox{ as }\eps\to 0\hskip.1cm.$$
It remains to prove the result in $B_0\left(\delta\right)$. Thanks to (\ref{eq3.24}), one can easily check that 
\begin{equation}\label{eq3.26bis}
\frac{u_\eps}{U_\eps} \to 1 +\sqrt{\lambda} \frac{G_2}{G_1} \hbox{ in }L^\infty\bigl(B_0\left(\delta\right)\bigr)\hbox{ as }\eps\to 0
\end{equation}
so that we can write that 
$$ 3u_\eps^5 = 3 U_\eps^5 + 15 U_\eps^4 \left(u_\eps-U_\eps\right) + 30 U_\eps^3\left(u_\eps-U_\eps\right)^2 + O\left(U_\eps^2 \left\vert u_\eps - U_\eps\right\vert^3\right)\hskip.1cm.$$
Using again (\ref{eq3.25}), we deduce that 
\begincal 
3 u_\eps^5 & = & 3 U_\eps^5 -15 \eps^{\frac{1}{2}} \gamma(0)U_\eps^4 + 15 \left(\lambda\eps\right)^{\frac{1}{2}} U_\eps^4 x^i\partial_i G_2(0)\\
&& + 15 \gamma(0) U_\eps^4 V_\eps + 15 U_\eps^4 \varphi_\eps \left(h(0)W_\eps + \gamma(0)^2 Y_\eps +Z_\eps\right)\\
&&+ 30 \gamma(0)^2 U_\eps^3\left(V_\eps-\eps^{\frac{1}{2}}\right)^2 + o\left(U_\eps\right)
\fincal
in $B_0\left(\delta\right)$. Thanks to (\ref{eq3.24}), to (\ref{eq3.25}) and to (\ref{eq3.26}), we can also write that 
\begincal 
\Delta u_\eps + h u_\eps &=& 3U_\eps^5 + 15 \gamma(0)U_\eps^4 V_\eps+ 15 \left(\lambda \eps\right)^{\frac{1}{2}}\varphi_\eps U_\eps^4 x^i \partial_i G_2(0)\\
&&- 15 \gamma(0)\eps^{\frac{1}{2}} U_\eps^4  + 15 U_\eps^4 \varphi_\eps \left(h(0)W_\eps + \gamma(0)^2 Y_\eps +Z_\eps\right)\\
&& + 30 \gamma(0)^2 \left(2\eps^{-1} G^{-4} U_\eps^7 - G^{-2} U_\eps^5\right) + h(0)\left(\eps^2 - 8 G^{-2}\right) U_\eps^5 \\
&& + h(0) \varphi_\eps   \left(8\eps U_\eps^3 -9\eps^2 U_\eps^5\right)+  30\gamma(0)^2 \varphi_\eps\left(\eps^3 U_\eps^7 - \eps^4 U_\eps^9\right)\\
&& +18 \left(\varphi_\eps-1\right) U_\eps^5 \vert x\vert \partial_i\gamma(0)x^i+o\left(U_\eps\right)\hskip.1cm.
\fincal 
Combining these two last equations, we get that 
\begincal 
\Delta u_\eps + h u_\eps -3u_\eps^5 &=& 
- 30 \gamma(0)^2 \left(U_\eps^3\left(V_\eps-\eps^{\frac{1}{2}}\right)^2-\varphi_\eps \eps^3 U_\eps^7 \right.\\
&&\quad \left. +\varphi_\eps \eps^4 U_\eps^9
-2\eps^{-1} G^{-4} U_\eps^7 + G^{-2} U_\eps^5\right)\\
&& + h(0)\left(\eps^2U_\eps^5 - 8 G^{-2}U_\eps^5+ 8\varphi_\eps \eps U_\eps^3 - 9 \varphi_\eps \eps^2 U_\eps^5\right)+ o\left(U_\eps\right)\\
&&+ 18\left(\varphi_\eps-1\right) \vert x\vert x^i\partial_i\gamma(0)U_\eps^5  
\fincal
It remains to remark using (\ref{eq3.23bis}) that 
\begincal 
&&U_\eps^3\left(V_\eps-\eps^{\frac{1}{2}}\right)^2-\varphi_\eps \eps^3 U_\eps^7 +\varphi_\eps \eps^4 U_\eps^9
-2\eps^{-1} G^{-4} U_\eps^7 + G^{-2} U_\eps^5\\
&&\quad =   \eps^2 G^{-2} U_\eps^9\left(1-\varphi_\eps\right)\\
&&\quad = -\frac{\eps^2}{\ln \eps^2} \ln \left(1+\eps^{-2} G^{-2}\right) G^{-2} U_\eps^9\\
&&\quad = -\frac{U_\eps}{\ln \eps^2} \eps^6 \ln \left(1+\eps^{-2} G^{-2}\right) G^{-2}\left(\eps^2+G^{-2}\right)^{-4} \\
&&\quad = O\left(\frac{U_\eps}{\ln \frac{1}{\eps}}\right)=o\left(U_\eps\right)\hskip.1cm,
\fincal
that 
\begincal 
&&\eps^2U_\eps^5 - 8 G^{-2}U_\eps^5+ 8\varphi_\eps \eps U_\eps^3 - 9 \varphi_\eps \eps^2 U_\eps^5\\
&&\quad = \left(-9\eps^2 U_\eps^5 + 8 \eps U_\eps^3\right)\left(\varphi_\eps-1\right)\\
&&\quad = \frac{U_\eps}{\ln \eps^2} \ln \left(1 +\eps^{-2}G^{-2}\right) \left(-9 \left(1+ \eps^{-2}G^{-2}\right)^{-2} + 8 \left(1 + \eps^2 G^{-2}\right)^{-1}\right)\\
&&\quad  = O\left(\frac{U_\eps}{\ln \frac{1}{\eps}}\right)=o\left(U_\eps\right)
\fincal
and that 
\begincal 
\left(\varphi_\eps-1\right)G^{-1} x^i\partial_i\gamma(0)U_\eps^5 &=& 
O\left( \frac{U_\eps}{\ln \eps^2}\ln \left(1 +\eps^{-2}G^{-2}\right) \eps^{-2} G^{-2} \left(1 + \eps^{-2}G^{-2}\right)^{-2}\right)\\
&=&O\left(\frac{U_\eps}{\ln \frac{1}{\eps}}\right)=o\left(U_\eps\right)
\fincal
to conclude thanks to (\ref{eq3.26bis}) that (\ref{eq3.11}) holds in $B_0\left(\delta\right)$ for this choice of $u_\eps$. As already said, this proves that (\ref{eq3.11}) holds for $u_\eps$ given by (\ref{eq3.14}) and this ends the proof of the theorem. \hfill $\blacksquare$

\medskip As already said, this result implies theorem \ref{thminstability}.

\section{Appendix}

\subsection{A general simple lemma on functions}\label{generallemma}

We prove a new version of the simple Lemma 1.1 of \cite{DruetHebey08}, replacing the compact manifold $M$ by a domain $\Omega$ in $\R^n$. 

\begin{lemma}
\label{iso0} 
Let $\Omega$ be a smooth bounded domain of $\rn$. Let $u\in C^1\left(\Omega\right)$ be a function positive in the interior and null on the boundary. Assume that 
$$\{x\in \Omega \hbox{ s.t. } \nabla u (x)=0 \hbox{ and } d(x,\partial\Omega) u^2(x)\geq 1 \}\neq\emptyset \hskip.1cm.$$
Then there exist $N \in \N^{*}$ and $N$ critical points of $u$, denoted by $(x_{1},\dots, x_{N})$, such that 
\be
\begin{split}
&d(x_{i},\partial\Omega ) u(x_{i})^2  \geq 1 \hbox{ for all } i\in [1,N],\\
&\vert x_{i}-x_{j}\vert u(x_{i})^2  \geq 1 \hbox{ for all } i \not= j \in  [1,N] 
\end{split}
\ee
and 
\be
\left( \min_{i\in [1,N]} \vert x_{i} -x \vert  \right) u(x)^2 \leq 1
\ee
for all critical points $x$ of $u$ such that  $d(x,\partial\Omega) u(x)^2 \geq 1$.
\end{lemma}

\medskip {\bf Proof of Lemma \ref{iso0}.} Let $\mathcal{C}_u$ be the set of critical points of $u$. Thanks to the Hopf Lemma, it is clear that $\mathcal{C}_u$ is a compact set of $\Omega$.  We let
$$K_0=\{x\in\mathcal{C}_u \hbox{ s.t. } d(x,\partial\Omega)u^2(x)\geq 1 \}$$
and we assume that $K_0\neq \emptyset$. We let $x_1\in K_0 $ and $K_1 \subset K_0$ be such that 
\be
u(x_1)= \max_{K_0} u
\ee 
and
\be 
K_1=\left\{ x\in K_0 \hbox{ s.t. } \vert x_1 -x\vert u(x)^2 \geq 1 \right\} \hskip.1cm.
\ee
Then we proceed by induction. Assuming we have constructed $K_0 \supset \dots \supset K_p$ and $x_1, \dots, x_p$ such that $x_i \in K_{i-1}$ for all $i\in [1,p]$, we let $x_{p+1}\in K_p$ and $K_{p+1} \subset K_{p}$ be such that 
\be
u(x_{p+1})= \max_{K_p} u
\ee 
and
\be 
K_{p+1}=\left\{ x\in K_p \hbox{ s.t. } \vert x_{p+1} -x\vert u(x_{p+1})^2 \geq 1\hbox{ and } \min_{i\in [1,p]} \vert x - x_i \vert u(x)^2 \geq 1 \right\} \hskip.1cm.
\ee
We claim that, at some step in the process, $K_p = \emptyset$. In order to prove it, we remark that at each step of the construction,
\beq
\label{iso01}
\vert x_i - x_j \vert u (x_i)^2 \geq 1 \hbox{ for all } i\not=j \in [1,p],
\eeq  
which will prove the claim, since $\Omega$ is bounded. We prove (\ref{iso01}) by induction. Let $p\geq1$. By definition, for all $x\in K_p$, we have
 \be
 \vert x_i -x \vert u(x)^2 \geq 1\hbox{ for all } i\in[1,p] \hskip.1cm.
 \ee
This holds in particular for $x=x_{p+1}$. Then,  for all $x\in K_p$, we also easily check that
 \be
 \vert x_i -x \vert u(x_i)^2 \geq 1\hbox{ for all } i\in[1,p] \hskip.1cm,
 \ee 
which is also true for $x_{p+1}$, and proves (\ref{iso01}). Let $N\in \N^{*}$ be such that $K_N=\emptyset$. We claim that
\beq
\label{iso04}
\left( \min_{i\in [1,N]} \vert x_{i} -x \vert  \right) u^2(x) \leq 1
\eeq
for all $x\in K_0$ which,  together with (\ref{iso01}), will end the proof of  the lemma. Let $x\in K_0$.  Since $K_N=\emptyset$, there exists $p$ such that $x\in  K_{p-1}$ and $x\not\in K_{p}$. Then  we have either
\be
 \vert x_{p} -x\vert u^2(x_{p}) < 1
\ee
or
\be
 \min_{i\in [1,p]} \vert x - x_i \vert u(x) <  1\hskip.1cm.
\ee
In the second case, (\ref{iso04}) is clearly true while in the first, using the definition of $x_p$, we have that
\be
 \vert x_{p} -x\vert u^2(x) \leq \vert x_{p} -x\vert u^2(x_{p}) < 1\hskip.1cm,
\ee
which proves that (\ref{iso04}) also holds. As already sais, this proves the lemma.\hfill$\blacksquare$

\subsection{ Green function of $\Delta  + h$}\label{greenfct}

We prove here some basic estimates on Green's functions of operators $\Delta+h$ where $h$ is of low regularity.  

\begin{lemma}\label{estimg}
Let $\Omega$ be a smooth bounded domain of $\R^3$. Let $h \in L^p\left(\Omega\right)$ for some $p>3$. Then there exists $\delta>0$ such that if
\beq
\Vert h_{-}\Vert_{\frac{3}{2}} < \delta \hskip.1cm,
\eeq
then the operator $\Delta + h$ admits a positive Green function $\mathcal{G}_{h}$ which verifies the following estimates~:  
\be
\left\vert \vert x-y\vert \mathcal{G}_{h}(x,y) -\frac{1}{\omega_2} \right\vert \leq C \vert x-y \vert 
\ee 
and
 \be
 \left\vert \vert x-y\vert^2 \vert\nabla \mathcal{G}_{h}(x,y)\vert -\frac{1}{\omega_2} \right\vert \leq C \vert x-y \vert 
 \ee 
for all $x\not=y \in \Omega$,  where C is a positive constant depending only on $\Omega$, $\Vert h \Vert_p$ and $\delta$.
\end{lemma}

\medskip {\bf Proof of lemma \ref{estimg}.} We divide the proof into three steps.

\medskip {\it Step 1 :  $\Delta + h$ is coercive if $\Vert h_{-}\Vert_{\frac{3}{2}}$ is small enough -} Let $u\in H_{0}^{1}(\Omega)$. We write that 
\be
\io \left(\vert \nabla u \vert^{2} +h u^2\right)\, dx \geq  \io \left(\vert \nabla u \vert^{2} -h_{-} u^2\right)\, dx
\geq \Vert \nabla u \Vert^{2}_{2} -\Vert h_{-}\Vert_{\frac{3}{2}} \Vert u\Vert_{6}
\ee
thanks to H\"older's inequalities. One can then use Sobolev's embeddings and the fact that $\Vert h_{-}\Vert_{\frac{3}{2}}$ is small to conclude this first step.

\medskip {\it Step 2 : Existence and a priori estimate -} Let $\mathcal{G}(x,y)$ be the Green function of the Laplacian. Then solving
\be
\begin{split}
\Delta_{y} \mathcal{G}_{h}(x,y) + h \mathcal{G}_{h}(x,y) &= \delta_{x} \hbox{ in } \Omega,\\
\mathcal{G}_{h}(x,y) &= 0 \hbox{ on } \partial\Omega,
\end{split}
\ee
is equivalent to solving
\be
\begin{split}
\Delta_{y} \beta (x,y) + h \beta(x,y) &= -h\mathcal{G}(x,y),\\
\beta(x,y)&=0 \hbox{ on } \partial\Omega.
\end{split}
\ee
Since $h\in L^p\left(\Omega\right)$ for some $p>3$, there exists $q>\frac{3}{2}$ such that $h\mathcal{G}(x,.)\in L^q\left(\Omega\right)$. The existence of $\beta$ follows from the coercivity of $\Delta +h $ and  Lax-Milgram theorem. Moreover, using again the coercivity of $\Delta+h$ and Sobolev's embeddings, we get that 
\be
\frac{1}{C} \Vert \nabla \beta \Vert_2^2 \leq \int_{\Omega} (\vert\nabla\beta\vert^2 + h \beta^2 )dx = \int_{\Omega} -h\mathcal{G} \beta dx \leq \Vert h \mathcal{G}\Vert_{\frac{3}{2}} \Vert \beta \Vert_3 \leq C \Vert \nabla \beta \Vert_{2},
\ee
for some  $C>1$ depending only on $\Vert h \Vert_p$, $\Vert h_{-} \Vert_\frac{3}{2}$ and $\Omega$. This  gives an a priori bound on $\Vert \nabla\beta \Vert_2 $.

\medskip {\it Step 3 : estimates and positivity -} Thanks to the previous Step, there exists $C>0$ which depends only on $\Vert h \Vert_p$ and $\Vert h_{-} \Vert_\frac{3}{2}$, and $q>\frac{3}{2}$ such that
\be
\Vert h(\beta + G(x,.)) \Vert_{q} \leq C\hskip.1cm.
\ee
Now, thanks to standard elliptic theory (see for instance theorem 9.13 of \cite{GilbargTrudinger}), we see that $\beta\in L^{\infty}$ and 
\be
\Vert \beta\Vert_{\infty} \leq C
\ee
where $C$ is a positive constant which depends only on $\Vert h \Vert_p$ and $\Vert h_{-} \Vert_\frac{3}{2}$. This proves the first estimate of the lemma. The second follows by standard elliptic theory. Positivity of the Green function is an easy consequence of the coercivity of the operator $\Delta+h$. \hfill$\blacksquare$

\subsection{General Poho\v{z}aev's identities}\label{Pohozaevidentities}

For the sake of completeness, we derive here several forms of the classical Poho\v{z}aev identity \cite{Pohozaev} we used in this paper. Assume that $u $ is a $C^2$- solution of
 \be
\Delta u = u^5 -hu \hbox{ in } \Omega \hskip.1cm.
\ee
Multiplying this equation by $<x,\nabla u>$ and integrating by parts, one easily gets that 
\beq
 \label{pohozaev1}
 \frac{1}{2} \io \left( h u^2 +  h<x,\nabla u^2> \right) dx = 
 B_1 +B_2 ,
 \eeq
 where 
 \be
 \begin{split}
 B_1&=  \ibo \left( <x, \nabla u> \partial_\nu u  + \frac{1}{2} u\partial_\nu u- <x,\nu> \frac{\vert \nabla u\vert^2}{2} \right)d\sigma \hbox{ and }\\
B_2 &= \ibo <x,\nu>   \frac{u^6}{6}  d\sigma \hskip.1cm.
 \end{split}
 \ee
Hence, if $u=0$ on $\partial\Omega$, we get that 
\beq
\label{pohozaev2}
\frac{1}{2} \io h\left( u^2 + <x,  \nabla u^2 > \right) dx = 
\ibo  <x,\nu> \left(\partial_\nu u \right)^{2} d\sigma \hskip.1cm.
\eeq
Integrating by parts again, we get the Poho\v{z}aev identity in its usual form~: 
\beq
\label{pohozaev3} 
\io \left( h + \frac{<x,\nabla h>}{2} \right)u^2 dx = -
\ibo  <x,\nu> \left(\partial_\nu u \right)^{2} d\sigma \hskip.1cm.
\eeq 
In a similar way, multiplying the equation by $\nabla u$ and integrating by parts, one can derive the following Poho\v{z}aev's identity~:
\beq
\label{pohozaev4}
\ibo \left( \frac{\vert \nabla u \vert^2}{2} \nu  - \partial_\nu u \nabla u  +   \frac{u^6}{6}\nu \right) d\sigma  = \io h \frac{\nabla u^2 }{2} dx\hskip.1cm.
\eeq

\subsection{ Poho\v{z}aev's identity for Green functions}\label{pohog}

In this section, we prove a useful Poho\v{z}aev identity for a sum of Green's functions. First of all, we easily derive the following Lemma from standard elliptic theory~:
 
\begin{lemma}
\label{decomp}
Let $\Omega$ be a smooth bounded domain in $\R^3$. Let $y\in \Omega$ and let $g$ be a weak solution in $H^1\left(\Omega\right)$ of 
\be
\Delta g + h g = -\frac{h}{\omega_{2}\vert x-y\vert} \hbox{ in } \Omega\hskip.1cm.
\ee
Then $g$ is continuous and can be written as 
\beq
g(x)= g(y) + \frac{h(y)}{2} \vert x-y\vert + \gamma_{y}(x) \hbox{ in }\Omega
\eeq
where $\gamma_{y} \in C^1(\Omega)$ satisfies $\gamma_{y}(y)=0$.
\end{lemma}

\medskip Applying the previous decomposition lemma to Green's functions, we get the following Poho\v{z}aev identity on the regular parts of them.

\begin{lemma}
\label{pohoglemma} 
Let $\Omega$ be a smooth bounded domain in $\R^3$, star-shaped with respect to $0$ and let $h\in C^1(\Omega)$ which satisfies (\ref{mainhyp}). Let $\mathcal{G}_{h}$ be the Green function of $\Delta +h$. Let also $N\in \N^{*}$, $x_{1},\dots, x_{N}\in \Omega$, $\lambda_{1}, \dots, \lambda_{N}$ some positive real numbers and
\be
G(x)=\sum_{i=1}^{N} \lambda_{i} \mathcal{G}_{h}(x,x_i)\hskip.1cm.
\ee
Then, using lemma \ref{decomp}, we write $G$ in a neighbourhood of $x_i$ as 
\be
G(x)=\frac{\lambda_i}{\omega_{2}\vert x - x_i \vert} + m_i +\lambda_i \frac{h(x_i)}{2} \vert x-x_{i}\vert +\gamma_{i}(x)
\ee
where $m_i\in \R$ and $\gamma_{i}\in C^1(\Omega)$ satisfies $\gamma_i(0)=0$. Then we have that 
\be
\sum_{i=1}^{N} \lambda_{i} \left(m_i + 2<x_i, \nabla \gamma_{x_{i}}(x_{i})>\right) <0 \hskip.1cm.
\ee
\end{lemma}

\medskip{\bf Proof of lemma \ref{pohog}.} We let $\delta>0$ be such that the $B(x_{i}, \delta)$ are disjoint and do not intersect the boundary of $\Omega$ and we set 
$$\Omega_\delta = \Omega\setminus \left\{ \bigcup_{i=1}^{N}B(x_{i},\delta) \right\}\hskip.1cm.$$ 
Multiplying the equation satisfied by $G$ by $<x,\nabla G>$ and after some integrations by parts, we obtain that 
\be
\begin{split}
& \int_{\Omega_\delta} \left(\frac{1}{2}<x,\nabla h> + h\right) G^2 \, dx\\
&\quad = \int_{\partial\Om} \left(\frac{1}{2}<x,\nu>\left(\left\vert \nabla G\right\vert^2 + hG^2\right) - \left(\left(x,\nabla G\right)+\frac{1}{2} G\right)\partial_\nu G \right)\, d\sigma\\
&\qquad - \sum_{i=1}^N\int_{\partial B\left(x_i,\delta\right)}  \left(\frac{1}{2}<x,\nu>\left(\left\vert \nabla G\right\vert^2 + hG^2\right) - \left(\left(x,\nabla G\right)+\frac{1}{2} G\right)\partial_\nu G \right)\, d\sigma
\end{split}
\ee
where $\nu$ denotes the outer normal to $\partial\Om$ and to $\partial B\left(x_i,\delta\right)$ respectively. Noting that $G=0$ on $\partial\Om$, we have that 
\be
\begin{split}
&\int_{\partial\Om} \left(\frac{1}{2}<x,\nu>\left(\left\vert \nabla G\right\vert^2 + hG^2\right) - \left(\left(x,\nabla G\right)+\frac{1}{2} G\right)\partial_\nu G \right)\, d\sigma\\
&\quad = -\frac{1}{2}\int_{\partial\Om} <x,\nu>\left\vert \nabla G\right\vert^2\, d\sigma <0
\end{split}
\ee
since $\Omega$ is star-shaped. Since $h$ satisfies (\ref{mainhyp}), we arrive to 
\be 
\sum_{i=1}^N\int_{\partial B\left(x_i,\delta\right)}  \left(\frac{1}{2}<x,\nu>\left(\left\vert \nabla G\right\vert^2 + hG^2\right) - \left(\left(x,\nabla G\right)+\frac{1}{2} G\right)\partial_\nu G \right)\, d\sigma \le -C_0
\ee
where $C_0$ is independent of $\delta$. It is easily checked that 
\be
\int_{\partial B\left(x_i,\delta\right)} <x,\nu>hG^2\, d\sigma\to 0\hbox{ as }\delta\to 0\hskip.1cm.
\ee
In order to estimate the remaining terms, we write that 
\be 
\begin{split}
&\int_{\partial B\left(x_i,\delta\right)}  \left(\frac{1}{2}<x,\nu>\left\vert \nabla G\right\vert^2  - \left(\left(x,\nabla G\right)+\frac{1}{2} G\right)\partial_\nu G \right)\, d\sigma\\
&\quad = \int_{\partial B\left(x_i,\delta\right)}  \left(\frac{1}{2}<x-x_i,\nu>\left\vert \nabla G\right\vert^2  - \left(\left(x-x_i,\nabla G\right)+\frac{1}{2} G\right)\partial_\nu G \right)\, d\sigma\\
&\qquad + \int_{\partial B\left(x_i,\delta\right)}  \left(\frac{1}{2}<x_i,\nu>\left\vert \nabla G\right\vert^2  - \left(x_i,\nabla G\right)\partial_\nu G \right)\, d\sigma \hskip.1cm.
\end{split}
\ee
Then, thanks to the expansion of $G$ in a neighbourhood of $x_i$, one can easily check that
\be
\int_{\partial B\left(x_i,\delta\right)}  \left(\frac{1}{2}<x-x_i,\nu>\left\vert \nabla G\right\vert^2  - \left(\left(x-x_i,\nabla G\right)+\frac{1}{2} G\right)\partial_\nu G \right)\, d\sigma \to \frac{\lambda_i}{2}m_i
\ee
and 
\be 
\int_{\partial B\left(x_i,\delta\right)}  \left(\frac{1}{2}<x_i,\nu>\left\vert \nabla G\right\vert^2  - \left(x_i,\nabla G\right)\partial_\nu G \right)\, d\sigma \to \lambda_i <x_i,\nabla \gamma_{x_i}\left(x_i\right)>
\ee
as $\delta\to 0$. Combining the above results gives the desired inequality.\hfill$\blacksquare$  

\bibliographystyle{amsplain}
\bibliography{DruetLaurainJEMS}

\end{document}